\documentclass[fleqn,final]{zzz}
\usepackage{amssymb}
\usepackage{amsmath, amsfonts}
\usepackage{mathrsfs}
\usepackage{color}
\usepackage{tikz}
\usepackage{float}
\usepackage{graphicx}
\usepackage{rotating}
\usepackage[pdftex]{hyperref}
\usepackage{enumerate}
\usepackage{soul}
\usepackage{comment}
\usepackage{xparse}
\usepackage{chngcntr}
\usepackage{apptools}
\usepackage{caption}
\usepackage[normalem]{ulem}
\usepackage{subcaption}
\usepackage{arydshln}
\theoremstyle{definition}
\newtheorem{ex}{Example}

\numberwithin{proposition}{section}
\numberwithin{remark}{section}
\numberwithin{equation}{section}

\DeclareMathOperator{\li}{Li}

\DeclareMathOperator{\arctanh}{arctanh}
\newcommand{\rmd}{\mathrm}

\newcommand{\G}{\textbf{\textup{G}}}

\newcommand{\img}{\mathrm{i}}
\newcommand{\B}{\mathrm{B}}
\newcommand{\eul}{\mathrm{e}}
\hypersetup{
colorlinks=true,
urlcolor=blue,
linkcolor=red,
citecolor=green
}

\NewDocumentCommand{\qfrac}{smm}{%
 \dfrac{\IfBooleanT{#1}{\vphantom{\big|}}#2}{\mathstrut #3}%
}

\makeatletter
\def\eqnarray{\stepcounter{equation}\let\@currentlabel=\theequation
\global\@eqnswtrue
\tabskip\@centering\let\\=\@eqncr
$$\halign to \displaywidth\bgroup\hfil\global\@eqcnt\z@
  $\displaystyle\tabskip\z@{##}$&\global\@eqcnt\@ne
  \hfil$\displaystyle{{}##{}}$\hfil
  &\global\@eqcnt\tw@ $\displaystyle{##}$\hfil
  \tabskip\@centering&\llap{##}\tabskip\z@\cr}

\def\endeqnarray{\@@eqncr\egroup
      \global\advance\c@equation\m@ne$$\global\@ignoretrue}

\def\@yeqncr{\@ifnextchar [{\@xeqncr}{\@xeqncr[5pt]}}
\makeatother

\newcommand{\C}{{\mathbb C}}
\newcommand{\N}{{\mathbb N}}

\newcommand{\R}{{\mathbb R}}

\definecolor{darkgreen}{rgb}{0.0, 0.21, 0.06}

\begin{document}

\renewcommand{\PaperNumber}{***}

\FirstPageHeading

\ShortArticleName{%
New closed forms for a class of digamma series and integrals
}

\ArticleName{%
New closed forms for a class of digamma series and integrals
}

\Author{Abdulhafeez A. Abdulsalam\,$^\dag\!\!\ $}

\AuthorNameForHeading{A.~A.~Abdulsalam}
\Address{$^\dag$ Department of Mathematics, University of Ibadan, Ibadan, Oyo, Nigeria}
\EmailD{hafeez147258369@gmail.com} 


\ArticleDates{Received \today~in final form ????; 
Published online ????}

\Abstract{The pursuit of closed forms for infinite series has long been a focal point of research. In this paper, we contribute to this endeavor by presenting closed forms for the class of digamma series:
\[\sum_{k=1}^\infty \frac{\psi\left(\frac{2k+2n+5}{4}\right) - \psi\left(\frac{2k+2n+3}{4}\right)}{(2k + \alpha)^2}, \quad \sum_{k=1}^\infty (-1)^k \frac{\psi\left(\frac{2k+2n+5}{4}\right) - \psi\left(\frac{2k+2n+3}{4}\right)}{(2k + \alpha)^2},\]
for all non-negative integers $\alpha$ and $n$. In addition to providing closed forms for these series, we unveil new identities for various generalized digamma series in the elegant form $a_0 + a_1 \pi + a_2 \pi^2 + a_3 \pi^3$, where $a_0, \ldots, a_3$ are real-valued constants determined by our formulas. Furthermore, we present ten definite integrals over the interval $(0, 1)$ that have not been previously studied in the literature and appear to be nearly impossible to evaluate. Combining these series and integrals can lead to the discovery of even more new results. Our findings contribute to the study of closed forms for infinite series and integrals, offering novel results and potential avenues for further exploration.}


\Classification{11M35, 33B15}
\allowdisplaybreaks
\section{Introduction}
In this paper, we delve into the study of infinite series of the forms
\begin{align}\label{eq1int}
\Theta_1(n, \alpha) &:= \sum_{k=1}^\infty \frac{f(k, n)}{(2k+\alpha)^2},\\
\label{eq2int}
\Theta_2(n, \alpha) &:= \sum_{k=1}^\infty (-1)^k \frac{f(k, n)}{(2k+\alpha)^2},
\end{align}
where $n$ and $\alpha$ are non-negative integers, and the function $f(k, n)$ is given by:
\begin{equation}\label{eq3int}
f(k, n) = \psi\left(\frac{2k+2n+5}{4}\right) - \psi\left(\frac{2k+2n+3}{4}\right).
\end{equation}
Notably, considering the series \eqref{eq1int} and \eqref{eq2int} with $2k$ in the denominator adds more strength to our investigation compared to considering them with $k$. The simplest instances of the series in \eqref{eq1int} and \eqref{eq2int} arise when $f(k, n) = 1$ and $m = 1$. In this particular case, the first series can be expressed in terms of the famous mathematical constant, $\pi$, yielding $\sum_{k=1}^\infty \frac{1}{(2k+1)^2} = \frac{\pi^2}{8} - 1$. On the other hand, the second series, $\sum_{k=1}^\infty \frac{(-1)^k}{(2k+1)^2}$, leads to $\G - 1$, where $\G = \sum_{k=0}^\infty \frac{(-1)^k}{(2k+1)^2}$ represents Catalan's constant ($\G \approx 0.9159655941$). The latter series cannot be expressed in terms of other mathematical constants and holds significance as a unique constant itself. Taking $f(k, n)$ as defined in \eqref{eq3int}, we observe intriguing similarities with our previous findings. Both $\Theta_1(n, \alpha)$ and $\Theta_2(n, \alpha)$ can be represented in terms of known mathematical constants. However, for the series $\Theta_2(n, \alpha)$, an additional constant emerges that defies further reduction to other known constants, namely $\Im(\li_3(1 + \img)) \approx 1.2670834419$, where $\Im(z)$ denotes the imaginary part of $z$, $\li_3(z)$ is the trilogarithm function, and $\img$ represents the complex unit, that is, $\img = \sqrt{-1}$.\\\\
The motivation for this research emerged from an attempt to evaluate one of V\u{a}lean's integral using a new approach. V\u{a}lean's first integral \cite[\S1.38, pp.~50, (1.176)]{bib9} is given by
\begin{equation}\label{Valeanint}
\begin{split}
2\int_0^1 \frac{\arctan{x}\ln(1+x)}{x}\, \rmd{d}x + \int_0^1 \frac{\arctan{x}\ln\left(1+x^2\right)}{x} \, \rmd{d}x \\
\quad + 4\int_0^1 \frac{x\arctan{x} \ln{x}}{1+x^2} \, \rmd{d}x = 4\G\ln{2} - \frac{\pi^3}{16}.
\end{split}
\end{equation}
During the process of integration by parts and utilizing the series expansion of the $\arctan$-function, the series $\Theta_2(n, \alpha)$ emerges for the specific values $n=0$ and $m=1$. This occurrence sparked our interest and prompted further exploration into generalizations of the forms $\Theta_1(n, \alpha)$ and $\Theta_2(n, \alpha)$, as defined earlier. In our pursuit, we employ a unique approach to evaluate V\u{a}lean's first integral and successfully prove the closed forms for the remaining integrals in his book \cite[\S1.38, pp.~50, (1.174), (1.175), (1.177)]{bib9}. Remarkably, our approach enables us to discover new closed forms for definite integrals over the interval $(0, 1)$ that are nearly impossible to evaluate.\\

In summary, the results established in this article are outlined as follows. In section \ref{section31}, we provide closed forms for $\Theta_1(n, 2m-1)$ and $\Theta_2(n, 2m-1)$, for all positive integers $m$. In sections \ref{sec3.6} and \ref{sec3.7}, we provide closed forms for $\Theta_1(n, 2m)$ and $\Theta_2(n, 2m)$, respectively for all non-negative integers $m$. In sections \ref{sec3.2} and \ref{sec3.4}, we prove the four V\u{a}lean integrals \cite[\S1.38, pp.~50, (1.174)--(1.177)]{bib9}. In section \ref{sec3.5}, we provide ten new definite integrals. In section \ref{sec3.8}, we provide new identities for some generalized digamma series.  The Computer Algebra System (CAS) software used to verify our results is \textsf{Mathematica~13}.

\section{Notations and Definitions}\label{sec2}
In this manuscript, we consistently use the natural logarithm denoted by $\ln{z}$, referring to the principal branch of the logarithm, where $\Im \ln{z} \in (-\pi, \pi]$. The principal branch is chosen since, for instance, for $|z| < 1$, $\ln(1-z)$ is multivalued, while its Taylor series $-\sum_{k=1}^\infty \frac{z^k}{k}$ is single-valued. To maintain clarity and conciseness, we adopt the following abbreviated notations: $ \overline{z}$ denotes the complex conjugate of $z$; $\Re(z)$ represents the real part of $z$; $\gamma \approx 0.5772156649$ denotes Euler's constant; $\mathrm{e} \approx 2.71828182845$ represents Euler's number; $H_n := \sum_{k=1}^n \frac{1}{k}$ represents the $n$-th harmonic number; $\mathbb{N}_0 := \mathbb{N} \cup \{0\}$ (where $n\mathbb{N}_0$ denotes the set of all elements of $\mathbb{N}_0$ multiplied by $n$); $\mathbb{Z}$, $\mathbb{Q}$, $\mathbb{R}$, and $\mathbb{C}$ represent the sets of integers, rational, real, and complex numbers, respectively.\\\\
We denote the gamma \cite[(5.2.1)]{bib23}, digamma \cite[(5.2.2)]{bib23}, trigamma, tetragamma, pentagamma, and hexagamma functions \cite[\S 5.15]{bib23} of argument $z$ as $\Gamma(z)$, $\psi(z)$, $\psi_1(z)$, $\psi_2(z)$, $\psi_3(z)$, and $\psi_4(z)$, respectively, where $\psi_n(z) := \frac{\rmd{d}^{n+1}}{\rmd{d}z^{n+1}}\ln{\Gamma(z)}$, and $n\in\mathbb{N}_0$. The function $\psi_n(z)$ is generally known as the polygamma function, and $\psi(z) := \psi_0(z)$.\\\\
The digamma function has the series representation \cite[\S1.7(6)]{bib2}
\[\psi(z) = -\gamma +\sum_{k=0}^\infty \left(\frac{1}{k+1} - \frac{1}{k+z}\right), \quad z \in \mathbb{C}\setminus- \mathbb{N}_0.\]
When $z = k+1$ (a positive integer), the digamma function becomes \cite[\S1.7.1(9)]{bib2}
\[\psi(k+1) = -\gamma + H_k, \quad k \in \mathbb{N}.\]
The recurrence relation for the digamma function is given by \cite[(5.5.2)]{bib23}
\begin{equation}\label{recc}
\psi(z+1) = \psi(z) + \frac{1}{z}.
\end{equation}
The reflection formula for the polygamma function is given by \cite[(5.15.6)]{bib23}
\begin{equation}\label{polyrefl}
\psi_n(z) + (-1)^{n+1} \psi_n(1-z) = -\pi\frac{\rmd{d}^n}{\rmd{d}z^n}\left(\cot(\pi z)\right), \quad n \in \mathbb{N}_0,\,z \in \C \setminus-\mathbb{N}_0.
\end{equation}
The duplication formula for the $\psi(z)$ is \cite[(5.5.8)]{bib23}
\begin{equation}\label{dupl}
\psi\left(z + \frac{1}{2}\right) = 2\psi(2z) - \psi(z) - \ln{4}, \quad z \in \C \setminus-\mathbb{N}_0.
\end{equation}
Euler's beta function is defined as \cite[(5.12.1)]{bib23}:
\[\B(a, b) := \int_0^1 t^{a-1} (1-t)^{b-1} \, \rmd{d}t = \frac{\Gamma(a)\Gamma(b)}{\Gamma(a+b)}, \quad a, b \in \C \setminus-\mathbb{N}_0.\]
The Lerch transcendent is defined as \cite[(24.14.1)]{bib23}:
$$\Phi(z, s, a) := \sum_{n=0}^{\infty} \frac{z^n}{(n + a)^s}, \quad \lvert z\rvert \leq 1, \Re\,s > 1, a \not\in -\mathbb{N}_0.$$
The polygamma function has the series representation $\psi_n(z) = (-1)^{n-1} n! \Phi(1, n+1, z)$, where $n \in \N$. The Dirichlet eta function is defined as $\eta(n) := \Phi(-1, n, 1)$ \cite[\S1.12(2)]{bib2}, where $\Re n > 0$. The Riemann zeta function \cite[\S 25.2]{bib23} and the Hurwitz zeta function \cite[\S 25.11]{bib23} are, respectively, defined as:
$$\zeta(s) := \sum_{n=1}^{\infty} \frac{1}{n^s}, \quad \zeta(s, z) := \sum_{n=0}^{\infty} \frac{1}{(n+z)^s},$$
where $z \not\in -\mathbb{N}_0$, $\Re\,s > 1$. The domain $\Re\,s > 1$ can be extended to $s \in \mathbb{C}\setminus \{1\}$ through analytic continuation, using for instance, the Hermite integral representation for the Hurwitz zeta function \cite[(25.11.29)]{bib23}. The relationship between the Dirichlet eta function and the Riemann zeta function is given by $\eta(n) = \left(1- 2^{1-n}\right)\zeta(n)$ \cite[\S1.12(2)]{bib2}. The polylogarithm function \cite[\S 25.12(ii), \S 25.14.3]{bib23}, $\li_s(z)$, is defined as: $\li_s(z) := z\Phi(z, s, 1)$, where $\Re s > 1,\, |z| \leq 1$. The dilogarithm function, $\li_2(z)$, has the integral representation \cite[(25.12.2)]{bib23}
\begin{equation}\label{intrep}
\li_2(z) = -\int_0^z \frac{\ln(1-t)}{t}\, \rmd{d}t, \quad z \in \C \setminus(1, \infty).
\end{equation}
The domain $|z| \leq 1$ can be extended to the entire complex plane through analytic continuation. This can be achieved using, for instance, the integral representation \cite[(25.14.6)]{bib23}
\[\li_s(z) = \frac{1}{2}z + \int_0^\infty \frac{z^{t+1}}{(1+t)^s} \, \rmd{d}t - 2z\int_0^\infty \frac{\sin\left(t\ln{z} - s\arctan{t}\right)}{\left(1+t^2\right)^{\frac{s}{2}} \left(e^{2\pi t} - 1\right)} \, \rmd{d}t, \quad \Re(s) > 0\,\,\textup{if}\,\, z\in \C\setminus[1, \infty).\]

\section{Results}
\subsection{Closed forms for $\Theta_1(n, 2m-1)$ and $\Theta_2(n, 2m-1)$}\label{section31}
In this subsection, we provide closed for $\Theta_1(n, 2m-1)$ and $\Theta_2(n, 2m-1)$, where $m$ is a positive integer. We require the following lemmas.
\begin{lemma}[Euler's reflection formula] Let $z\in \C$, where $z\neq 0, 1$. Then \cite[(25.12.6)]{bib23}
\begin{equation}\label{lem1}
\li_2(z) + \li_2(1-z) = \frac{\pi^2}{6} - \ln{z}\ln(1-z).
\end{equation}
\end{lemma}

\begin{lemma}[Landen's identity] Let $z \in \C\setminus[1, \infty)$. Then \cite[(25.12.3)]{bib23}
\begin{equation}\label{lem2}
\li_2(z) + \li_2\left(\frac{z}{z-1}\right) = -\frac{1}{2}\ln^2(1- z).
\end{equation}
\end{lemma}

\begin{lemma}Let $z \in \C\setminus(-\infty, 0]$. Then 
\begin{equation}\label{lem3}
\int_0^z \frac{\ln{t}\ln(1-t)}{t}\, \rmd{d}t = \li_3(z) - \li_2(z)\ln{z}.
\end{equation}
\end{lemma}

\begin{lemma}Let $z\in \C\setminus[1, \infty)$. Then 
\begin{equation}\label{lem4}
\int_0^z \frac{\ln^2(1-t)}{t}\, \rmd{d}t = 2\zeta(3) - 2\li_3(1- z) + \ln(1- z)\left(\ln{z} \ln(1- z) + 2\li_2(1- z)\right).
\end{equation}
\end{lemma}

\begin{lemma}Let $z \in \C\setminus[1, \infty)$. Then 
\begin{equation}\label{lem5}
\li_3(z) + \li_3(1-z) + \li_3\left(\frac{z}{z-1}\right) = \frac{1}{6}\ln^2(1- z)\left(\ln(1-z) - 3\ln{z}\right) + \frac{\pi^2}{6}\ln(1-z) + \zeta(3).
\end{equation}
\end{lemma}

\begin{lemma}Let $z\in \C\setminus(-\infty,-1]$. Then 
\begin{equation}\label{lem5b}
\begin{split}
\int_0^z \frac{\ln^2(1+t)}{t}\, \rmd{d}t &= 2\zeta(3) - \frac{2}{3}\ln^3(1 + z) - 2\li_3\left(\frac{1}{1+z}\right) - 2\ln(1+z)\li_2\left(\frac{1}{1+z}\right)
\\&\quad+\ln{z}\ln^2(1+z).
\end{split}
\end{equation}
\end{lemma}

\begin{proof}
The proof of \eqref{lem1} follows from utilizing the transformation $z \longmapsto 1-z$ in the integral representation \eqref{intrep} of $\li_2(z)$. The proof of \eqref{lem2} can be derived by differentiating the integral representation of $\li_2\left(\frac{z}{z-1}\right)$ using \eqref{intrep}, then separating by partial fractions and integrating again. As for the proof of \eqref{lem3}, it can be obtained by expressing \eqref{lem3} as
\[\int_0^z \frac{\ln{t}\ln(1-t)}{t}\, \rmd{d}t  = \frac{\ln^2{z}\ln(1-z)}{2} + \frac{1}{2}\int_0^z \frac{\ln^2{t}}{1-t} \, \rmd{d}t,\]
and subsequently performing term-wise integration. The proof of \eqref{lem4} follows from expressing \eqref{lem4} as 
\[\int_0^z \frac{\ln^2(1-t)}{t}\, \rmd{d}t = \ln^2(1-z)\ln{z} + 2\zeta(3) -2 \int_0^{1-z} \frac{\ln{t}\ln(1-t)}{t}\, \rmd{d}t,\]
and then using \eqref{lem3}. Changing the variable $z$ in \eqref{lem2} to $t$, multiplying both sides of \eqref{lem2} by $\frac{1}{t(1-t)}$, and integrating from $t=0$ to $z$, we obtain
\[\li_3\left(\frac{z}{z-1}\right) + \li_3(z) + \int_0^z \frac{\li_2(t)}{1-t}\, \rmd{d}t = -\frac{1}{2}\int_0^z \frac{\ln^2(1-t)}{t}\, \rmd{d}t - \frac{1}{2}\int_0^z\frac{\ln^2(1-t)}{1-t}\, \rmd{d}t.\]
Integrating by parts, we have
\begin{equation}\label{eqlme3}
\li_3\left(\frac{z}{z-1}\right) + \li_3(z) - \ln(1-z)\li_2(z) = \frac{1}{2}\int_0^z \frac{\ln^2(1-t)}{t}\, \rmd{d}t + \frac{1}{6}\ln^3(1-z).
\end{equation}
Using \eqref{lem3} and \eqref{lem1} in \eqref{eqlme3}, we conclude the proof of \eqref{lem5}. By changing the variable $z$ in \eqref{lem5b} to $t$, making the transformation $t \longmapsto t - 1$, and then further transforming $t \longmapsto \frac{1}{t}$, the integral \eqref{lem5b} takes the form
\begin{equation}\label{term1s}
\int_0^z \frac{\ln^2(1+t)}{t}\, \rmd{d}t = \int_{\frac{1}{1+z}}^1 \frac{\ln^2{t}}{t}\,\rmd{d}t + \int_{\frac{1}{1+z}}^1 \frac{\ln^2{t}}{1-t}\,\rmd{d}t.
\end{equation}
The first integral in \eqref{term1s} simplifies to $\frac{1}{3} \ln^3(1+z)$, and the second integral is evaluated using term-wise integration, leading to the final result.
\end{proof}

\begin{lemma} \label{lemma2} Let $z \in \C\setminus[1,\infty)$. Then
\begin{equation}\label{lem6}
\sum_{k=1}^\infty \frac{H_k}{(k+1)^2} z^{k+1} = \zeta(3) + \li_2(1-z)\ln(1-z) - \li_3(1-z) + \frac{1}{2}\ln{z}\ln^2(1-z),
\end{equation}
\begin{equation}\label{lem7}
\sum_{k=1}^\infty \frac{\psi_1(k+1)}{k+1}z^{k+1} = -z\zeta(2) + \li_3(z) + \int_0^1 \frac{\ln{t}\ln(1-zt)}{1-t} \, \rmd{d}t.
\end{equation}
\end{lemma}

\begin{proof}
Using the integral representation for the harmonic numbers \cite[\S1.7.2(13)]{bib2}
\begin{equation}\label{harmonic1}
H_k = \int_0^1 \frac{1-t^k}{1-t} \, \rmd{d}t,
\end{equation}
we have that
\begin{equation}\label{eqpl1}
\begin{split}
\sum_{k=1}^\infty \frac{H_k}{k+1} z^{k+1} &= \int_0^1 \frac{\ln(1-zt)}{t} + \frac{\ln\left(\frac{1-zt}{1-z}\right)}{1-t} \, \rmd{d}t = -\li_2(z) - z \int_0^1 \frac{\ln(1-t)}{1-zt}\,\rmd{d}t
\\&=  -\li_2(z) + \int_0^{\frac{z}{z-1}} \frac{\ln\left(\frac{z-1}{z}t\right)}{1-t} \, \rmd{d}t = -\li_2(z) + \int_0^{\frac{z}{z-1}} \frac{\ln(1-t)}{t} \, \rmd{d}t 
\\&=  -\li_2(z) - \li_2\left(\frac{z}{z-1}\right) = \frac{1}{2}\ln^2(1-z).
\end{split}
\end{equation}
To complete the proof of \eqref{lem6}, we begin by changing the variable in \eqref{eqpl1} from $z$ to $t$, then divide both sides by $t$, integrate both sides from $t=0$ to $t=z$, and finally apply \eqref{lem4}.
\end{proof}
By expressing $H_k$ as $\gamma + \psi(k+1)$ and differentiating both sides with respect to $k$ in \eqref{harmonic1}, we obtain
\begin{equation}\label{psi1zz}
\sum_{k=1}^\infty \frac{\psi_1(k+1)}{k+1}z^{k+1} = z\int_0^1 \frac{\ln{t}}{1-t} \, \rmd{d}t + \int_0^1 \frac{\ln{t}\ln(1-zt)}{t(1-t)}\, \rmd{d}t.
\end{equation}
The proof of \eqref{lem7} is completed by splitting the second integral in \eqref{psi1zz}.
\begin{remark}If we substitute $z=1$ in \eqref{lem6} and make the transformation $t \longmapsto 1-t$ in the integral, we derive
\begin{equation}
\sum_{k=1}^\infty \frac{\psi_1(k+1)}{k+1} = 2\zeta(3) - \zeta(2).
\end{equation}
Substituting $z = \frac12$ in \eqref{lem1} and \eqref{lem4}, we have
\begin{align}
\li_2\left(\frac12\right) &= \frac{\pi^2}{12} - \frac{\ln^2{2}}{2} \label{dilog12},\\
\li_3\left(\frac12\right) &= \frac{\ln^3{2}}{6} - \frac{\pi^2}{12} \ln{2} + \frac{7}{8}\zeta(3).\label{trilog12}
\end{align}
The identities \eqref{dilog12} and \eqref{trilog12} are known in the literature. 
\end{remark}
The following lemma will be useful in the subsequent remark, and will also play an important role in proving Theorem \ref{thm2} and \eqref{Valeanint}.
\begin{lemma}Let $s \in \R$ and $z \in \C$. Then $\Re \li_s(z) = \Re \li_s( \overline{z})$ and $\Im \li_s(z) = -\Im \li_s( \overline{z})$.\label{lemconj1}
\end{lemma}
\begin{proof}
Let $s \in \R$, $z = x+ \img y\in \C$, where $x, y \in \R$, and $|z| \leq 1$. Then
\begin{equation}\label{polylogconj}
\begin{split}
\li_s(z)  &= \sum_{k=1}^\infty \frac{\left(x^2 + y^2\right)^{\frac{k}{2}} \left(\cos\left(k\arctan\frac{y}{x}\right) + \img \sin\left(k\arctan\frac{y}{x}\right)\right)}{k^s}, \\
\li_s(\overline{z}) &= \sum_{k=1}^\infty \frac{\left(x^2 + y^2\right)^{\frac{k}{2}} \left(\cos\left(k\arctan\frac{y}{x}\right) - \img \sin\left(k\arctan\frac{y}{x}\right)\right)}{k^s}.
\end{split}
\end{equation}
From \eqref{polylogconj}, we can deduce that $\Re \li_s(z) = \Re \li_s(\overline{z})$ and $\Im \li_s(z) = -\Im \li_s(\overline{z})$. These equations remain valid for all $z \in \mathbb{C}$ and $s \in \mathbb{R}$ due to analytic continuation. However, it is crucial to note that the equations are not valid if $s \in \mathbb{C}$ and $\Im(s) \neq 0$ (or $s\in \C \setminus\R$). To demonstrate this, consider letting $s = a + \mathrm{i}b$, where $a, b \in \mathbb{R}$ and $b \neq 0$. Then, the following holds
\begin{equation}\label{exception1}
\begin{split}
\li_s(z) &= \sum_{k=1}^\infty \frac{\left(x^2 + y^2\right)^{\frac{k}{2}} \left(\cos(b\ln{k}) - \img \sin(b\ln{k})\right) \left(\cos\left(k\arctan\frac{y}{x}\right) + \img \sin\left(k\arctan\frac{y}{x}\right)\right)}{k^a}, \\
\li_s(\overline{z}) &= \sum_{k=1}^\infty \frac{\left(x^2 + y^2\right)^{\frac{k}{2}} \left(\cos(b\ln{k}) - \img \sin(b\ln{k})\right) \left(\cos\left(k\arctan\frac{y}{x}\right) - \img \sin\left(k\arctan\frac{y}{x}\right)\right)}{k^a}.
\end{split}
\end{equation}
It is evident from \eqref{exception1} that $\Re \li_s(z) \neq \Re \li_s( \overline{z})$ and $\Im \li_s(z) \neq -\Im \li_s(\overline{z})$. This inequality holds true for all $z \in \mathbb{C}$ due to analytic continuation.
\end{proof}

\begin{remark}We require the following special values
\begin{equation}\label{li2i}
\li_2(\img) = \sum_{k=1}^\infty \frac{(-1)^k}{(2k)^2} - \img \sum_{k=1}^\infty \frac{(-1)^k}{(2k-1)^2} = -\frac{\pi^2}{48} + \img \G,\quad \li_2(-\img) = \frac{\pi^2}{48} - \img \G,
\end{equation}
\begin{equation}\label{li3i}
\li_3(\img) = \sum_{k=1}^\infty \frac{(-1)^k}{(2k)^3} - \img \sum_{k=1}^\infty \frac{(-1)^k}{(2k-1)^3} = -\frac{3}{32}\zeta(3) + \frac{\img \pi^3}{32},\quad \li_3(-\img)  = -\frac{3}{32}\zeta(3) - \frac{\img \pi^3}{32}.
\end{equation}
It is worth noting that the evaluation of the second series in \eqref{li3i} involves transforming the series into its tetragamma form and then substituting $n=2$ and $z=\frac{1}{4}$ into \eqref{polyrefl}. By employing \eqref{lem1} and \eqref{li2i}, we obtain
\begin{equation}\label{remark2}
\li_2(1 - \img) = \frac{\pi^2}{16} - \img\left(\frac{\pi}{4}\ln{2} + \G\right), \quad \li_2(1 + \img) = \frac{\pi^2}{16} + \img\left(\frac{\pi}{4}\ln{2} + \G\right).
\end{equation}
Substituting $z=1-\img$ in \eqref{lem5} and using \eqref{li3i}, we derive
\begin{equation}\label{lemli3ref}
\li_3(1- \img) + \li_3(1+\img) = \frac{\pi^2}{16} \ln{2} + \frac{35}{32}\zeta(3).
\end{equation}
\end{remark}
\newpage\noindent
In the following theorem, we provide the generalization $\Theta_1(n, 1)$.
\begin{theorem}Let $n \in \N$. Then \label{thm1}
\begin{align*}
&\sum_{k=1}^\infty \frac{\psi\left(\frac{2k+2n+5}{4}\right) - \psi\left(\frac{2k+2n+3}{4}\right)}{(2k+1)^2}
\\&\quad= \begin{cases} 
\begin{aligned}[b]
- 4 + (1 - \G)\pi  + \frac{\pi^2}{4}\ln{2} + \frac{7}{4}\zeta(3),
\end{aligned} & \textup{if $n=0$,} \\ \ \\
\begin{aligned}[b]
\frac{5}{3} - (1 - \G)\pi  + \left(1-\ln{2}\right)\frac{\pi^2}{4} - \frac{7}{4}\zeta(3), 
\end{aligned} & \textup{if $n=1$},\\ \ \\
\begin{aligned}[b]
 - \frac{14}{5} + (1 - \G)\pi + \left(\ln{2} - \frac{1}{2}\right)\frac{\pi^2}{4} + \frac{7}{4}\zeta(3), 
\end{aligned} & \textup{if $n=2$},\\ \ \\
\begin{aligned}[b]&3 -\frac{1}{8}\sum_{j=1}^{\frac{n-1}{2}} \frac{(4j+1)\left(H_{2j-1} - 2 H_{4j-1}\right)  }{j^2 (2j+1)^2}
\\&- \sum_{j=1}^{\frac{n-1}{2}} \frac{1}{(2j+1)^2(4j+1)}+ 4 \sum_{j=1}^{n}   \frac{(-1)^j}{2j+1} - (1 - \G)\pi  
\\&+ \left(1- \ln{2} - \frac{1}{2}\sum_{j=1}^{\frac{n-1}{2}} \frac{1}{j(2j+1)}\right)\frac{\pi^2}{4} - \frac{7}{4}\zeta(3),
\end{aligned} & \begin{tabular}{cc} \textup{if $n$ is odd}\\ \textup{and $n \neq1$,}\end{tabular}\\ \ \\
\begin{aligned}[b]
&- \frac{10}{3} - \frac{1}{8}\sum_{j=1}^{\frac{n-2}{2}} \frac{(3+4j)\left(H_{2j-1} - 2 H_{4j-1}\right)}{(j+1)^2(2j+1)^2}
\\&-  \sum_{j=1}^{\frac{n-2}{2}} \frac{4j^3 + j^2 - 4j - 2}{(j+1)^2(2j+1)^2(4j+3)(4j+1)} - 4\sum_{j=1}^{n} \frac{(-1)^j}{2j+1}
\\&+ (1 - \G)\pi + \left(\ln{2} - \frac{1}{2}  - \frac{1}{2}\sum_{j=1}^{\frac{n-2}{2}} \frac{1}{(j+1)(2j+1)}\right) \frac{\pi^2}{4}  + \frac{7}{4}\zeta(3),
\end{aligned} & \begin{tabular}{cc} \textup{if $n$ is even}\\ \textup{and $n \neq2$.}\end{tabular}
\end{cases}
\end{align*}
\end{theorem}

\begin{proof}
We start by proving the theorem for the case $n=0$. Using the recurrence relation \eqref{recc} and the duplication formula \eqref{dupl}, we unveil the following identity
\begin{equation}\label{greresultar}
\psi\left(\frac{2k+5}{4}\right) - \psi\left(\frac{2k+3}{4}\right) = (-1)^{k-1}\pi + 4(-1)^k\sum_{j=0}^k \frac{(-1)^j}{2j+1}, \quad k \in \N.
\end{equation}
The identity \eqref{greresultar} compels us to begin the proof by considering the series
$$\sum_{k=1}^\infty \frac{(-1)^k}{(2k+1)^2}\sum_{j=1}^k \frac{(-1)^{j-1}}{2j+1}.$$
Interchanging the order of summation, we have
\begin{equation}\label{nr1}
\sum_{k=1}^\infty \sum_{j=1}^k \frac{(-1)^{k+j-1}}{(2k+1)^2(2j+1)} = \sum_{j=1}^\infty \sum_{k=j}^\infty \frac{(-1)^{k+j-1}}{(2k+1)^2(2j+1)}.
\end{equation}
Interchanging the roles of the dummy variables $j$ and $k$, we obtain
\begin{equation}\label{nr2}
\begin{split}
\sum_{k=1}^\infty \sum_{j=1}^k \frac{(-1)^{k+j-1}}{(2k+1)^2(2j+1)} &= \sum_{j=1}^\infty \sum_{k=1}^j \frac{(-1)^{k+j-1}}{(2j+1)^2 (2k+1)} 
\\& = \left(\frac{\pi}{4} - 1\right)(1- \G) + 1 - \frac{7}{8}\zeta(3) + \sum_{j=1}^\infty \sum_{k=j}^\infty \frac{(-1)^{k+j}}{(2j+1)^2 (2k+1)}.
\end{split}
\end{equation}
Adding \eqref{nr1} and \eqref{nr2}, it follows that
\begin{equation}\label{eqsp1}
\begin{split}
\sum_{k=1}^\infty \sum_{j=1}^k \frac{(-1)^{k+j-1}}{(2k+1)^2(2j+1)}  &= \frac{\pi}{8} - \frac{\pi\G}{8} + \frac{\G}{2} - \frac{7\zeta(3)}{16} + \sum_{j=1}^\infty \sum_{k=1}^\infty\frac{k  (-1)^{k}}{(2k + 2j +1)^2 (2j+1)^2}.
\end{split}
\end{equation}
Letting $S$ represent the final sum on the right-hand side of \eqref{eqsp1}, we have
\begin{equation}\label{nr3}
\begin{split}
S &= \sum_{j=1}^\infty \sum_{k=1}^\infty \left(\frac{(-1)^k}{4k(2j+1)^2} - \frac{(-1)^k}{2k(2j+1)(2j+2k+1)} +  \frac{(-1)^k}{4k(2k+2j+1)^2}\right)
\\&= \frac{-\pi^2}{32}\ln{2} - \frac{\ln{2}}{4} - \frac{\pi}{4} + 1  -  \frac{1}{2}\int_0^1 \sum_{k=1}^\infty \frac{(-1)^k }{k} x^{2k}  \sum_{j=0}^\infty \frac{x^{2j+1}}{2j+1} \, \frac{\rmd{d}x}{x} 
\\&\qquad+  \sum_{j=1}^\infty \sum_{k=1}^\infty\frac{(-1)^k}{4k(2k+2j+1)^2}
\\&= \frac{-\pi^2}{32}\ln{2} - \frac{\ln{2}}{4} - \frac{\pi}{4} + 1 + \frac{1}{2} \int_0^1 \frac{\ln\left(1+ x^2\right)}{x} \arctanh{x} \, \rmd{d}x 
\\&\qquad+ \sum_{j=1}^\infty \sum_{k=1}^\infty\frac{(-1)^k}{4k(2k+2j+1)^2}.
\end{split}
\end{equation}
Letting $S_1$ represent the last sum on the right-hand side of $S$, we have
\begin{equation}\label{nr4}
\begin{split}
S_1 &= \sum_{k=1}^\infty \sum_{j=k}^\infty \frac{(-1)^k}{4k(2j + 1)} \int_0^1 x^{2j}\, \rmd{d}x -  \sum_{k=1}^\infty \frac{(-1)^k}{4k(2k+1)^2}
\\&= \sum_{k=1}^\infty \sum_{j=0}^\infty\frac{(-1)^k}{4k(2j + 1)^2}  - \sum_{k=1}^\infty \sum_{j=0}^{k-1} \frac{(-1)^k}{4k(2j + 1)} \int_0^1 x^{2j}\, \rmd{d}x - \sum_{k=1}^\infty \frac{(-1)^k}{4k(2k+1)^2}
\\&= -\frac{\pi^2}{32}\ln{2} + \frac{\ln{2}}{4} + \frac{\G}{2} - 1 + \frac{\pi}{8} - \sum_{k=1}^\infty \sum_{j=0}^{k-1} \frac{(-1)^k}{4k(2j + 1)} \int_0^1 x^{2j}\, \rmd{d}x.
\end{split}
\end{equation}
Letting $S_2$ denote the final sum on the right-hand side of $S_1$, we have
\begin{equation}\label{nr5}
\begin{split}
S_2 &=  \sum_{k=1}^\infty \frac{(-1)^k}{4k} \int_0^1 \sum_{j=0}^{k-1} \frac{x^{2j+1}}{2j+1} \frac{\rmd{d}x}{x} = \sum_{k=1}^\infty \frac{(-1)^k}{4k} \int_0^1 \int_0^x \frac{1- y^{2k}}{1-y^2} \, \rmd{d}y \frac{\rmd{d}x}{x} 
\\&= \frac{1}{4}\int_0^1 \int_0^x \frac{-\ln{2} + \ln\left(1+y^2\right)}{1-y^2}\, \rmd{d}y \frac{\rmd{d}x}{x}  = \frac{1}{4}\int_0^1 \int_y^1 \frac{-\ln{2} + \ln\left(1+y^2\right)}{1-y^2}\, \frac{\rmd{d}x}{x}\, \rmd{d}y 
\\&= \frac{1}{4}\int_0^1 \frac{\ln{2} \ln{y} - \ln\left(1+y^2\right)\ln{y}}{1-y^2}\, \rmd{d}y.
\end{split}
\end{equation}
From \eqref{nr3}, \eqref{nr4}, \eqref{nr5}, it follows that
\begin{equation}\label{splitssum}
\begin{split}
S &= -\frac{\pi^2}{16}\ln{2} - \frac{\pi}{8} + \frac{\G}{2} + \frac{1}{2} \int_0^1 \frac{\ln\left(1+ x^2\right)}{x} \arctanh{x} \, \rmd{d}x  - \frac{\ln{2}}{4} \int_0^1 \frac{\ln{x}}{1 - x^2} \, \rmd{d}x 
\\&\quad+ \frac{1}{4}\int_0^1 \frac{\ln\left(1+ x^2\right)\ln{x}}{1-x^2}\, \rmd{d}x.
\end{split}
\end{equation}
If we let $I$ represent all the integrals in $S$, then
\begin{equation}\label{splitsint}
\begin{split}
I &= \frac{1}{4}\int_0^1 \frac{\ln\left(1 + x^2\right) \ln\left(1-x^2\right)}{x} \, \rmd{d}x - \frac{1}{4}\int_0^1 \frac{\ln\left(1 - x\right) \ln\left(1+x^2\right)}{x} \, \rmd{d}x 
\\&\quad+ \frac{1}{4}\int_0^1 \ln\left(1 + x^2\right) \left(\frac{\ln{x}}{1 - x^2} - \frac{\ln(1- x)}{x}\right) \rmd{d}x - \frac{\ln{2}}{4}\int_0^1 \frac{\ln{x}}{1-x^2} \rmd{d}x.
\end{split}
\end{equation}
For the third integral in \eqref{splitsint}, denoted by $I_3$, we have
\begin{align*}
I_3 &= \frac14\int_0^1  \ln\left(1 + x^2\right) \rmd{d}\left(\frac{1}{4}\li_2\left(x^2\right) + \ln{x}\arctanh{x}\right) 
\\&= \frac12\int_0^1 \frac{x \ln{x} \ln(1-x)}{1+x^2} \, \rmd{d}x - \frac14\int_0^1 \frac{x \ln{x} \ln\left(1-x^2\right)}{1+x^2} \, \rmd{d}x - \frac{1}{16}\int_0^1 \frac{\ln(1-x)\ln(1+x)}{x} \, \rmd{d}x.
\end{align*}
This implies
\begin{equation}\label{eqhi}
\begin{split}
I &= \frac{1}{16}\int_0^1 \frac{\ln(1-x)\ln(1+x)}{x} \, \rmd{d}x - \frac{1}{4}\int_0^1 \frac{\ln(1- x)\ln\left(1+x^2\right)}{x} \, \rmd{d}x 
\\&\quad+ \frac{1}{2}\int_0^1 \frac{x \ln{x} \ln(1-x)}{1+x^2} \, \rmd{d}x - \frac{1}{4}\int_0^1 \frac{x \ln{x} \ln\left(1-x^2\right)}{1+x^2} \, \rmd{d}x - \frac{\ln{2}}{4}\int_0^1 \frac{\ln{x}}{1-x^2} \rmd{d}x.
\end{split}
\end{equation}
Taking limits as $z\to1$ \eqref{lem4}, we obtain
\[\int_0^1 \frac{\ln^2\left(1-x^2\right)}{x} \, \rmd{d}x = \frac{1}{2}\int_0^1 \frac{\ln^2(1- x)}{x}\, \rmd{d}x = \zeta(3), \quad \int_0^1 \frac{\ln^2(1- x)}{x}\, \rmd{d}x  = 2\zeta(3).\]
Taking limits as $z\to1$ in \eqref{lem5b} and using \eqref{dilog12}, \eqref{trilog12}, we have
\[\int_0^1 \frac{\ln^2(1+x)}{x}\, \rmd{d}x = \frac{1}{4}\zeta(3).\]
Therefore, we have
\begin{equation}\label{hisol1}
\begin{split}
\int_0^1 \frac{\ln(1-x)\ln(1+x)}{x}\, \rmd{d}x &= \frac{1}{2}\int_0^1 \frac{\left(\ln^2\left(1-x^2\right) - \ln^2(1 - x) - \ln^2(1+x)\right)}{x} \, \rmd{d}x
\\&= \frac{1}{2}\left(\zeta(3) - 2\zeta(3) - \frac{\zeta(3)}{4}\right) = -\frac{5}{8}\zeta(3).
\end{split}
\end{equation}
The integral \eqref{hisol1} is also evaluated in \cite[pp.~73--74]{bib32}. For the third integral in \eqref{eqhi}, we have
\begin{equation}\label{hisol3}
\begin{split}
\int_0^1 \frac{x \ln{x} \ln\left(1-x\right)}{1+x^2} \, \rmd{d}x &= \zeta(3) + \int_0^1 \int_0^1\left(\frac{\alpha^2 \ln{x}}{(1+\alpha^2)(1 - \alpha x)} + \frac{(1+\alpha x)\ln{x}}{(1+\alpha^2)(1+x^2)}\right) \rmd{d}x\, \rmd{d}\alpha
\\&= \zeta(3) - \frac{3\pi^2}{32} \ln{2} - \frac{\pi \G}{4} - \frac{1}{2}\int_0^1 \frac{\ln(1-x)\ln(1 + x^2)}{x}\, \rmd{d}x.
\end{split}
\end{equation}
For the fourth integral in \eqref{eqhi}, we have
\begin{equation}\label{hisol4}
\begin{split}
\int_0^1 \frac{x \ln{x} \ln\left(1-x^2\right)}{1+x^2} \, \rmd{d}x &= -\frac{1}{4}\int_0^1 \int_0^1 \left(\frac{\ln{x}}{(1-\alpha x)(1 + \alpha)}  -  \frac{\ln{x}}{(1+x)(1+\alpha)} \right) \rmd{d}x \, \rmd{d}\alpha 
\\&= \frac{1}{4}\zeta(3) - \frac{\pi^2}{16}\ln{2} - \frac{1}{4}\int_0^1 \frac{\ln(1-x)\ln(1+x)}{x} \, \rmd{d}x
\\&= \frac{13}{32}\zeta(3) - \frac{\pi^2}{16}\ln{2}.
\end{split}
\end{equation}
For the last integral in \eqref{eqhi}, we have
\begin{equation}\label{hisol5}
\begin{split}
\int_0^1 \frac{\ln{x}}{1-x^2} \, \rmd{d}x = -\sum_{k=0}^\infty \frac{1}{(2k+1)^2} = -\frac{3}{4}\zeta(2) = -\frac{\pi^2}{8}.
\end{split}
\end{equation}
Substituting \eqref{hisol1}, \eqref{hisol3}, \eqref{hisol4}, and \eqref{hisol5} into \eqref{eqhi}, and subsequently using the result in \eqref{splitssum}, we obtain
\begin{equation}\label{penssum}
S = -\frac{\pi^2}{16}\ln{2} - \frac{\pi}{8} - \frac{\pi \G}{8} + \frac{\G}{2} + \frac{23}{64}\zeta(3) - \frac{1}{2}\int_0^1 \frac{\ln(1-x)\ln(1 + x^2)}{x}\, \rmd{d}x.
\end{equation}\noindent
To evaluate the integral in \eqref{penssum}, we will employ the use of two complex valued integrals. Firstly, writing the series expansion for $1/(1-\img x)$ and subsequently using the beta functions in integrating term-wise, we obtain
\begin{equation}\label{frstint}
\begin{split}
\int_0^1 \frac{\ln{x}\ln(1-x)}{1-\img x} \, \rmd{d}x &= \sum_{k=0}^\infty \img^k \left(\frac{H_{k+1}}{(k+1)^2} - \frac{\psi_1(k+1)}{k+1} + \frac{1}{(k+1)^3}\right)
\\&= \sum_{k=1}^\infty \frac{\img^k H_k}{(k+1)^2} - 2\img \li_3(\img) - \sum_{k=0}^\infty \frac{\img^k \psi_1(k+1)}{k+1}.
\end{split}
\end{equation}
Using Lemma \ref{lemma2} and \eqref{remark2} in \eqref{frstint}, we obtain
\begin{equation}\label{impint}
\begin{split}
\int_0^1 \frac{\ln{x}\ln(1-x)}{1-\img x} \, \rmd{d}x &= -\frac{\G}{2}\ln{2} - \frac{\pi}{16}\ln^2{2} + \img \left(\frac{\G \pi}{4} - \frac{\pi^2}{32}\ln{2} + \li_3(1-\img) - \frac{29}{32}\zeta(3)\right)
\\&\quad+\int_0^1 \frac{\ln{x}\ln(1-x)}{1-\img x} \, \rmd{d}x+ \img\int_0^1 \frac{\ln(1-x)\ln\left(1+x^2\right)}{x}\, \rmd{d}x.
\end{split}
\end{equation}
Secondly, by replacing $1 - \img x$ in \eqref{frstint} with its conjugate and performing analogous operations, we obtain
\begin{equation}\label{impint2}
\begin{split}
\int_0^1 \frac{\ln{x}\ln(1-x)}{1+\img x} \, \rmd{d}x &= -\frac{\G}{2}\ln{2} - \frac{\pi}{16}\ln^2{2} + \img \left(-\frac{\G \pi}{4} + \frac{\pi^2}{32}\ln{2} - \li_3(1+\img) + \frac{29}{32}\zeta(3)\right)
\\&\quad+\int_0^1 \frac{\ln{x}\ln(1-x)}{1+\img x} \, \rmd{d}x- \img\int_0^1 \frac{\ln(1-x)\ln\left(1+x^2\right)}{x}\, \rmd{d}x.
\end{split}
\end{equation}
Utilizing \eqref{impint}, \eqref{impint2}, and applying \eqref{lemli3ref}, we deduce
\begin{equation}\label{lnpnexsq}
\begin{split}
\int_0^1 \frac{\ln(1-x)\ln\left(1+x^2\right)}{x}\, \rmd{d}x &= -\frac{\G \pi}{2} + \frac{\pi}{16}\ln^2{2} - \left(\li_3(1-\img) + \li_3(1+\img)\right) + \frac{29}{16}\zeta(3)
\\& = -\frac{\G \pi}{2} + \frac{23}{32}\zeta(3).
\end{split}
\end{equation}
By substituting \eqref{lnpnexsq} into \eqref{penssum} and subsequently using the result in \eqref{eqsp1}, we finally get
\begin{equation}\label{sumjkcat}
\sum_{k=1}^\infty \frac{(-1)^k}{(2k+1)^2}\sum_{j=1}^k \frac{(-1)^{j-1}}{2j+1} = \G - \frac{\pi^2}{16}\ln{2} - \frac{7}{16}\zeta(3).
\end{equation}
Expressing \eqref{sumjkcat} as
\begin{align*}
\sum_{k=1}^\infty \frac{(-1)^k}{(2k+1)^2}\sum_{j=1}^k \frac{(-1)^{j-1}}{2j+1} = \sum_{k=1}^\infty \frac{(-1)^k}{(2k+1)^2}\left(\sum_{j=1}^\infty \frac{(-1)^{j-1}}{2j+1} - \sum_{j=0}^\infty \frac{(-1)^{j+k}}{2j+2k+3}\right), 
\end{align*}
we deduce that
\begin{equation}\label{usres}
\sum_{k=1}^\infty  \frac{\psi\left(\frac{2k+5}{4}\right) - \psi\left(\frac{2k+3}{4}\right)}{(2k+1)^2} = \frac{\pi^2}{4}\ln{2} + \frac{7}{4}\zeta(3) + (1 - \G)\pi - 4.
\end{equation} 
Using the recurrence relation \eqref{recc} and the duplication formula \eqref{dupl}, we derive the following functional equations
\begin{equation}\label{reldig}
\begin{split}
&\psi\left(\frac{2k+2n+5}{4}\right) - \psi\left(\frac{2k+2n+3}{4}\right) 
\\&\qquad= \begin{cases}\begin{aligned}[b]&\frac{4}{2k+3} - \psi\left(\frac{2k+5}{4}\right) +\psi\left(\frac{2k+3}{4}\right)
\\&\qquad+ 4\sum_{j=1}^{\frac{n-1}{2}}\left(\frac{1}{2k+4j+3} - \frac{1}{2k+4j+1}\right),\end{aligned} & \textup{if $n$ is odd and $n \neq 1$},\\ \ \\
\begin{aligned}[b]&\psi\left(\frac{2k+5}{4}\right) - \psi\left(\frac{2k+3}{4}\right) 
\\&\qquad+ 4\sum_{j=0}^{\frac{n-2}{2}}\left(\frac{1}{2k+4j+5} - \frac{1}{2k+4j+3}\right),\end{aligned} & \textup{if $n$ is even}.
\end{cases}
\end{split}
\end{equation}
To conclude the proof of Theorem \ref{thm1}, we proceed as follows: For $n$ odd and $n\neq1$, and for $n$ even and $n\neq2$, we multiply both sides of \eqref{reldig} by $1/(2k+1)^2$, sum the resulting equation from $k=1$ to $\infty$, and then utilize \eqref{usres}. For $n=1$, we use \eqref{dupl}, and for $n=2$, we use \eqref{recc}, and then we apply \eqref{usres}.
\end{proof}
In the following theorem, we present the generalization $\Theta_2(n, 1)$.
\begin{theorem}Let $n \in \N$. Then \label{thm2}
\begin{align*}
&\sum_{k=1}^\infty (-1)^k \frac{\psi\left(\frac{2k+2n+5}{4}\right) - \psi\left(\frac{2k+2n+3}{4}\right)}{(2k+1)^2}
\\&\quad= \begin{cases}
\begin{aligned}[b]
-4 + \left(1- \frac{1}{4}\ln^2{2}\right)\pi - \frac{\pi^3}{8}  + 4\Im(\li_3(1+\img)), 
\end{aligned} & \textup{if $n=0$.}\\ \ \\ 
\begin{aligned}[b]
\frac{11}{3} - \left(3 - \frac{1}{2}\ln^2{2}\right)\frac{\pi}{2} + \frac{\pi^3}{8} + 2\G - 4\Im(\li_3(1+\img)) , 
\end{aligned} & \textup{if $n=1$}, \\ \ \\
\begin{aligned}[b]&5  - \sum_{j=1}^{\frac{n-1}{2}} \frac{1}{4j^2(4j+1)} + 4\sum_{j=1}^{n} \frac{(-1)^j}{2j+1} 
\\&+ \frac{1}{4} \sum_{k=1}^{\frac{n-1}{2}} \sum_{j=0}^{2k}  \frac{(-1)^j(8k^2+4k+1)}{k^2(2j+1)(2k+1)^2} 
\\&+ \left(2- \sum_{j=1}^{\frac{n-1}{2}} \frac{1}{j(2j+1)}\right)\G - \left(3 + \sum_{j=1}^{\frac{n-1}{2}} \frac{1}{(2j+1)^2} - \frac{1}{2}\ln^2{2}\right) \frac{\pi}{2}  
\\&+ \frac{\pi^3}{8} - 4\Im(\li_3(1+\img)),
\end{aligned} & \begin{tabular}{cc} \textup{if $n$ is odd}\\ \textup{and $n \neq1$.}\end{tabular},\\ \ \\
\begin{aligned}[b]&-4 + \frac{1}{4} \sum_{j=1}^{\frac{n}{2}} \frac{1}{j^2(4j-1)}  - 4\sum_{j=1}^{n} \frac{(-1)^j}{2j+1} 
\\&- \frac{1}{4} \sum_{k=1}^{\frac{n}{2}} \sum_{j=0}^{2k-2}  \frac{(-1)^j(8k^2-4k+1)}{k^2(2j+1)(2k-1)^2} 
\\&+ \left(1 + \frac{1}{2} \sum_{j=1}^{\frac{n}{2}} \frac{1}{(2j-1)^2} - \frac{1}{4}\ln^2{2}\right)\pi - \frac{\pi^3}{8} 
\\&- \G \sum_{j=1}^{\frac{n}{2}} \frac{1}{j(2j-1)} +  4\Im(\li_3(1+\img)),
\end{aligned} & \textup{if $n$ is even}.
\end{cases}
\end{align*}
\end{theorem}

\begin{proof}
By the series expansion of $\arctan{t}$, we have
\begin{equation}\label{eqfinalres1}
\begin{split}
\int_0^1 \frac{\arctan{t}\ln\left(1+ t^2\right)}{t} \rmd{d} t &= \sum_{k=0}^\infty \frac{(-1)^k}{2k+1}\left(\frac{\ln{2}}{2k+1} - \frac{2}{2k+1}\int_0^1 \frac{t^{2k+2}}{1+t^2}\, \rmd{d}t\right)
\\&= \G \ln{2} - \frac{1}{2}\sum_{k=0}^\infty (-1)^k \frac{\psi\left(\frac{2k+5}{4}\right) - \psi\left(\frac{2k+3}{4}\right)}{(2k+1)^2}.
\end{split}
\end{equation}
By the logarithmic representation of $\arctan{t}$, we have
\begin{equation}\label{arctanlogint}
\int_0^1 \frac{\arctan{t}\ln\left(1+ t^2\right)}{t} \rmd{d} t = \frac{\img}{2}\int_0^{\img} \frac{\ln^2(1-t)}{t} \, \rmd{d}t - \frac{\img}{2}\int_0^{-\img} \frac{\ln^2(1-t)}{t} \, \rmd{d}t.
\end{equation}
Using \eqref{lem4} and Lemma \ref{lemconj1}, \eqref{arctanlogint} becomes
\begin{equation}\label{finalres551}
\int_0^1 \frac{\arctan{t}\ln\left(1+ t^2\right)}{t} \rmd{d} t = \frac{\pi^3}{16} + \G \ln{2} + \frac{\pi}{8}\ln^2{2} - 2\Im(\li_3(1+\img)).
\end{equation}
Equating \eqref{eqfinalres1} and \eqref{finalres551}, a closed form for $\sum_{k=1}^\infty \frac{(-1)^k}{(2k+1)^2}\left(\psi\left(\frac{2k+5}{4}\right) - \psi\left(\frac{2k+3}{4}\right)\right)$ is established. The rest of the proof of other cases of Theorem \ref{thm2} is similar to the proof of corresponding cases in Theorem \ref{thm1}.
\end{proof}

\subsubsection{\large Generalization of the two dgamma series to two variables}
In the following theorem, we present a generalization of Theorems \ref{thm1} and \ref{thm2}, extending the results from $\Theta_1(n, 1)$ to $\Theta_1(n, 2m+1)$, and from $\Theta_2(n, 1)$ to $\Theta_2(n, 2m+1)$, respectively, where $m \in \mathbb{N}$. 
\begin{theorem}Let $m$ and $n$ be any two positive integers such that $n \geq m$. Then \label{thm3}
\begin{align*}
&\sum_{k=1}^\infty \frac{\psi\left(\frac{2k+2n+5}{4}\right) - \psi\left(\frac{2k+2n+3}{4}\right)}{(2k+2m+1)^2}
\\&\quad= \begin{cases} 
\begin{aligned}[b]
&-4  - 4 \sum_{k=1}^{m} \sum_{j=0}^k \frac{(-1)^{j+k}}{(2j+1)(2k+1)^2} +\left(1+\sum_{j=1}^{m} \frac{(-1)^j}{(2j+1)^2} - \G\right)\pi 
\\&+\frac{\pi^2}{4}\ln{2} + \frac{7}{4}\zeta(3),
\end{aligned} & \textup{if $n=m$}, \\ \ \\
\begin{aligned}[b]
&\frac{5}{3} + 4 \sum_{k=1}^{m} \sum_{j=0}^{k+1} \frac{(-1)^{j+k}}{(2j+1)(2k+1)^2} - \left(1 + \sum_{j=1}^{m} \frac{(-1)^j}{(2j+1)^2} - \G\right)\pi 
\\&+ \frac{\pi^2}{4}\left(1-\ln{2}\right) - \frac{7}{4}\zeta(3),
\end{aligned} & \textup{if $n=m+1$},\\ \\ \\
\begin{aligned}[b]
&\frac{-14}{5} - 4 \sum_{k=1}^{m} \sum_{j=0}^{k+2} \frac{(-1)^{j+k}}{(2j+1)(2k+1)^2} 
\\&+\left(1 + \sum_{j=1}^{m} \frac{(-1)^j}{(2j+1)^2} - \G\right)\pi + \left(\ln{2} - \frac{1}{2}\right)\frac{\pi^2}{4} + \frac{7}{4}\zeta(3),
\end{aligned}& \textup{if $n=m+2$}, \\ \ \\
\begin{aligned}[b]
&3  -\frac{1}{8}\sum_{j=1}^{\frac{n-m-1}{2}} \frac{(4j+1)\left(H_{2j-1} - 2 H_{4j-1}\right)  }{j^2 (2j+1)^2} 
\\&- \sum_{j=1}^{\frac{n-m-1}{2}} \frac{1}{(2j+1)^2(4j+1)}  - 4\sum_{k=1}^{m} \sum_{j=0}^{k+n-m} \frac{(-1)^{k+n-m + j}}{(2k+1)^2 (2j+1)} 
\\&+ 4 \sum_{j=1}^{n-m} \frac{(-1)^j}{2j+1} + \left(- 1 - \sum_{j=1}^{m} \frac{(-1)^{j+n-m-1}}{(2j+1)^2} + \G \right)\pi 
\\&+ \left(1- \ln{2} - \frac{1}{2}\sum_{j=1}^{\frac{n-m-1}{2}} \frac{1}{j(2j+1)}\right)\frac{\pi^2}{4} - \frac{7}{4}\zeta(3),
\end{aligned} & \begin{tabular}{cc}\textup{if $n-m$ is odd}\\\textup{and $n-m \neq 1$,}\end{tabular}\\ \ \\
\begin{aligned}[b]
&\frac{-10}{3}  -\frac{1}{8}\sum_{j=1}^{\frac{n-m-2}{2}} \frac{(4j+3)\left(H_{2j-1} - 2 H_{4j-1}\right)  }{(j+1)^2 (2j+1)^2} 
\\&- \sum_{j=1}^{\frac{n-m-2}{2}} \frac{4j^3 + j^2 - 4j - 2}{(j+1)^2(2j+1)^2(4j+3)(4j+1)} 
\\&- 4\sum_{k=1}^{m} \sum_{j=0}^{k+n-m} \frac{(-1)^{k+n-m + j}}{(2k+1)^2 (2j+1)} - 4\sum_{j=1}^{n-m} \frac{(-1)^j}{2j+1} 
\end{aligned}
\end{cases}
\end{align*}

\begin{align*}
\qquad\quad\begin{cases} 
\begin{aligned}[b]
&- \left(- 1 + \sum_{j=1}^{m} \frac{(-1)^{j+n-m-1}}{(2j+1)^2} + \G \right)\pi + \frac{\pi^2}{4}\left(\ln{2} - \frac12\right)  
\\&- \frac{\pi^2}{8}\sum_{j=1}^{\frac{n-m-2}{2}} \frac{1}{(j+1)(2j+1)} + \frac{7}{4}\zeta(3),
\end{aligned}  &  \begin{tabular}{cc}\textup{if $n-m$ is even}\\\textup{and $n-m \neq 2$.}\end{tabular}
\end{cases}
\end{align*}

\end{theorem}

\begin{proof}
By expressing the sum as
\begin{align*}
&\sum_{k=1}^\infty \frac{\psi\left(\frac{2k+2(n-m)+5}{4}\right) - \psi\left(\frac{2k+2(n-m)+3}{4}\right)}{(2k+1)^2} - \sum_{k=1}^{m}  \frac{\psi\left(\frac{2k+2(n-m)+5}{4}\right) - \psi\left(\frac{2k+2(n-m)+3}{4}\right)}{(2k+1)^2},
\end{align*}
and making use of \eqref{greresultar}, Theorem \ref{thm1}, and the proof methodology employed in proving Theorem \ref{thm1}, we conclude the proof.
\end{proof}

\begin{theorem}Let $m$ and $n$ be any two positive integers such that $n \geq m$. Then \label{thm4}
\begin{align*}
&\sum_{k=1}^\infty (-1)^k \frac{\psi\left(\frac{2k+2n+5}{4}\right) - \psi\left(\frac{2k+2n+3}{4}\right)}{(2k+2m+1)^2}
\\&\quad= \begin{cases}
\begin{aligned}[b]
&(-1)^m\left(-4 - 4 \sum_{k=1}^m \sum_{j=0}^k \frac{(-1)^j}{(2j+1)(2k+1)^2}  \right.
\\& \left. + \left(1 - \frac{\ln^2{2}}{4} + \sum_{k=1}^m \frac{1}{(2k+1)^2}\right)\pi  - \frac{\pi^3}{8} + 4\Im(\li_3(1+\img))\right), 
\end{aligned} & \textup{if $n=m$,}\\ \ \\ 
\begin{aligned}[b]
&(-1)^m \left(\frac{11}{3}  + 4\sum_{k=1}^m \sum_{j=0}^{k+1} \frac{(-1)^j}{(2j+1)(2k+1)^2} - \frac{3\pi}{2}  \right.
\\&\left.+\left(\frac{\ln^2{2}}{4} -\sum_{k=1}^m \frac{1}{(2k+1)^2}\right)\pi + \frac{\pi^3}{8} + 2\G - 4\Im(\li_3(1+\img))\right),
\end{aligned} & \textup{if $n=m+1$,}\\ \ \\
\begin{aligned}[b]&(-1)^m\left(5  - \sum_{j=1}^{\frac{n-m-1}{2}} \frac{1}{4j^2(4j+1)} - 4\sum_{k=1}^m \sum_{j=0}^{k+n-m} \frac{(-1)^{j+n-m}}{(2j+1)(2k+1)^2}   \right.
\\&+ \frac{1}{4} \sum_{k=1}^{\frac{n-m-1}{2}} \sum_{j=0}^{2k}  \frac{(-1)^j (8k^2+4k+1)}{k^2(2j+1)(2k+1)^2} 
\\&+ 4\sum_{j=1}^{n-m} \frac{(-1)^j}{2j+1} +\left(2- \sum_{j=1}^{\frac{n-m-1}{2}} \frac{1}{j(2j+1)}\right)\G 
\\&+\left(\frac{\ln^2{2}}{4} - \frac{3}{2} + (-1)^{n-m}\sum_{j=1}^m \frac{1}{(2j+1)^2} - \frac{1}{2}\sum_{j=1}^{\frac{n-m-1}{2}} \frac{1}{(2j+1)^2}\right) \pi
\\&\left.+ \frac{\pi^3}{8} - 4\Im(\li_3(1+\img))\right),
\end{aligned} & \begin{tabular}{cc}\textup{if $n-m$ is odd}\\\textup{and $n-m \neq 1$,}\end{tabular}
\end{cases}
\end{align*}

\newpage
\begin{align*}
\qquad\quad\begin{cases}
\begin{aligned}[b]&(-1)^m\left(-4 + \frac{1}{4} \sum_{j=1}^{\frac{n-m}{2}} \frac{1}{j^2(4j-1)} - 4\sum_{k=1}^m \sum_{j=0}^{k+n-m} \frac{(-1)^{j+n-m}}{(2j+1)(2k+1)^2}   \right.
\\&- \frac{1}{4} \sum_{k=1}^{\frac{n-m}{2}} \sum_{j=0}^{2k-2}  \frac{(-1)^j(8k^2-4k+1)}{k^2(2j+1)(2k-1)^2} - 4\sum_{j=1}^{n-m} \frac{(-1)^j}{2j+1}
\\&+ \left(1 +  (-1)^{n-m}\sum_{j=1}^m \frac{1}{(2j+1)^2}+ \frac{1}{2} \sum_{j=1}^{\frac{n-m}{2}} \frac{1}{(2j-1)^2}\right)\pi 
\\&\left.-  \frac{\pi}{4}\ln^2{2} - \frac{\pi^3}{8} - \G \sum_{j=1}^{\frac{n-m}{2}} \frac{1}{j(2j-1)} +  4\Im(\li_3(1+\img))\right),
\end{aligned} & \begin{tabular}{cc}\textup{if $n-m$}\\ \textup{ is even.}\end{tabular}
\end{cases}
\end{align*}
\end{theorem}

\begin{proof}
By expressing the sum as
\begin{align*}
&\sum_{k=1}^\infty (-1)^{k+m} \frac{\psi\left(\frac{2k+2(n-m)+5}{4}\right) - \psi\left(\frac{2k+2(n-m)+3}{4}\right)}{(2k+1)^2} 
\\&\qquad- \sum_{k=1}^{m}  (-1)^{k+m} \frac{\psi\left(\frac{2k+2(n-m)+5}{4}\right) - \psi\left(\frac{2k+2(n-m)+3}{4}\right)}{(2k+1)^2},
\end{align*}
and making use of \eqref{greresultar}, the rest of the proof is similar to the proof methodology of Theorem \ref{thm4}.
\end{proof}
The following identity is important in evaluating V\u{a}lean's integral.
\begin{lemma}Let $z \in \left(-\frac{\pi}{2}, \frac{\pi}{2}\right)$. Then \label{newlem1}
\begin{equation}\label{beauide}
\ln^2\left(2\cos{z}\right) = 2\sum_{k=1}^\infty \frac{(-1)^{k-1} \cos(2(k+1)z)}{k+1}H_k + z^2.
\end{equation}
\end{lemma}

\begin{proof}
By making use of $\eul^{\img 2z} + \eul^{-\img 2z} = 2\cos{(2z)}$ in $2\sum_{k=1}^\infty \frac{\cos(2(k+1)z)}{k+1}H_k$ and using \eqref{eqpl1}, we obtain
\begin{align*}
2\sum_{k=1}^\infty \frac{\cos(2(k+1)z)}{k+1}H_k &= \frac{1}{4}\ln^2\left(2 - 2\cos{(2z)}\right) - \arctan^2\left(\frac{\sin{(2z)}}{1-\cos{(2z)}}\right)
\\&= \frac{1}{4}\ln^2\left(2\sin^2{z}\right) - \arctan^2\left(\cot{z}\right)
\\&= \frac{1}{4}\ln^2\left(2\sin^2{z}\right)  - \left(\frac{\pi}{2} - z\right)^2.
\end{align*}
Replacing $z$ with $\frac{\pi}{2} - z$, we conclude the proof of \eqref{beauide}.
\end{proof}

\subsection{Proof of V\u{a}lean's first integral}\label{sec3.2}
In this subsection, we shall use Lemma \ref{newlem1} and Theorem \ref{thm2} to prove V\u{a}lean's first integral.
\begin{proof}
For brevity, denote V\u{a}lean's first integral in \eqref{Valeanint} as $\mathcal{V}_1$ and perform integration by parts on the third integral to obtain the expression
\begin{equation}\label{Valean2}
\mathcal{V}_1 = 2\int_0^1 \frac{\arctan{x}\ln(1+x)}{x}\, \rmd{d}x - \int_0^1 \frac{\arctan{x}\ln\left(1+x^2\right)}{x} \, \rmd{d}x  - 2\int_0^1 \frac{\ln\left(1+x^2\right) \ln{x}}{1+x^2} \, \rmd{d}x.
\end{equation}
By using the series expansion of $\arctan{x}$ as employed in \eqref{eqfinalres1}, we obtain the following results
\begin{equation}\label{arctansq1}
\begin{split}
\int_0^1 \frac{\arctan{x}\ln\left(1+ x^2\right)}{x} \rmd{d}x &= \G \ln{2} - 2 + \frac{\pi}{2} - \frac{1}{2}\sum_{k=1}^\infty (-1)^k \frac{\psi\left(\frac{2k+5}{4}\right) - \psi\left(\frac{2k+3}{4}\right)}{(2k+1)^2}
\\&= \G \ln{2} + \frac{\pi^3}{16} + \frac{\pi}{8}\ln^2{2} - 2\Im(\li_3(1+\img)),
\end{split}
\end{equation}
\begin{equation}\label{atanoneplx}
\int_0^1 \frac{\arctan{x}\ln\left(1+ x\right)}{x} \rmd{d}x = 2\G \ln{2} - 1 - \sum_{k=1}^\infty \frac{(-1)^k\left(H_{2k+1} - H_k\right)}{(2k+1)^2}. 
\end{equation}
With the substitution $x=\tan{\theta}$, we find the expression for the third integral in \eqref{Valean2}
\begin{equation}\label{lnln1}
\begin{split}
\int_0^1 \frac{\ln\left(1+x^2\right) \ln{x}}{1+x^2} \, \rmd{d}x &= \frac{1}{2}\int_0^{\frac{\pi}{2}} \ln^2(\sin{\theta}) \, \rmd{d}\theta + \ln{2}\int_0^{\frac{\pi}{2}} \ln(\sin{\theta}) \, \rmd{d}\theta - \frac{\pi}{4}\ln^2{2} 
\\&\quad+ 2\int_0^{\frac{\pi}{4}} \ln^2(\cos{\theta}) \, \rmd{d}\theta.
\end{split}
\end{equation}
By using beta functions, we derive the following results
\begin{align*}
\int_0^{\frac{\pi}{2}} \ln(\sin{\theta}) \, \rmd{d}\theta &= -\frac{\pi}{2}\ln{2}, \quad \int_0^{\frac{\pi}{2}} \ln^2(\sin{\theta}) \, \rmd{d}\theta = \frac{\pi}{2}\ln^2{2} + \frac{\pi^3}{24}.
\end{align*}
This leads to the result
\begin{align*}
\int_0^1 \frac{\ln\left(1+x^2\right) \ln{x}}{1+x^2} \, \rmd{d}x &= 2\int_0^{\frac{\pi}{4}} \ln^2(\cos{\theta}) \, \rmd{d}\theta - \frac{\pi}{2}\ln^2{2} + \frac{\pi^3}{48}.
\end{align*}
Using the following integrals
\begin{align*}
\int_0^{\frac{\pi}{4}} \ln\left(\sin{\theta}\cos{\theta}\right) \, \rmd{d}\theta &= -\frac{\pi}{2}\ln{2}, \quad \int_0^{\frac{\pi}{4}} \ln(\tan{\theta}) \, \rmd{d}\theta = -\G,
\end{align*}
it can be easily shown that
\begin{equation}\label{coscat1}
\int_0^{\frac{\pi}{4}} \ln(\cos{\theta}) \, \rmd{d}\theta = -\frac{\pi}{4}\ln{2} + \frac{\G}{2}.
\end{equation}
Integrating \eqref{beauide} from $z=0$ to $\frac{\pi}{4}$ and using \eqref{coscat1}, we obtain
\begin{equation}\label{cossq11s}
\int_0^{\frac{\pi}{4}} \ln^2(\cos{\theta}) \, \rmd{d}\theta = -\G\ln{2} + \frac{\pi}{4}\ln^2{2} + \frac{\pi^3}{192} + \sum_{k=1}^{\infty} \frac{(-1)^{k-1} H_{2k}}{(2k+1)^2}.
\end{equation}
Using this result in \eqref{lnln1}, we get
\begin{equation}\label{lnlnsq134}
\int_0^1 \frac{\ln\left(1+x^2\right) \ln{x}}{1+x^2} \, \rmd{d}x = -2\G \ln{2} + \frac{\pi^3}{32} + 2 \sum_{k=1}^{\infty} \frac{(-1)^{k-1} H_{2k}}{(2k+1)^2}.
\end{equation}
Using these results in \eqref{Valean2}, we finally obtain the expression for $\mathcal{V}_1$ as
\begin{equation}\label{vuleanf}
\begin{split}
\mathcal{V}_1 &= 7\G\ln{2} + 2 \sum_{k=1}^{\infty} \frac{(-1)^{k} H_{2k}}{(2k+1)^2} + 2\sum_{k=1}^{\infty} \frac{(-1)^k H_k}{(2k+1)^2} + 2 - 2\sum_{k=0}^\infty \frac{(-1)^{k}}{(2k+1)^3}
\\&\quad+ \frac{1}{2}\sum_{k=1}^\infty (-1)^k \frac{\psi\left(\frac{2k+5}{4}\right) - \psi\left(\frac{2k+3}{4}\right)}{(2k+1)^2} - \frac{\pi}{2} - \frac{\pi^3}{16} .
\end{split}
\end{equation}
Recall from \eqref{li3i} that
\begin{equation}\label{vul2}
\sum_{k=0}^\infty \frac{(-1)^{k}}{(2k+1)^3} = \frac{\pi^3}{32}.
\end{equation}
By writing $(-1)^k$ as $\cos(\pi k)$, and exploiting the fact that $\cos\left(\frac{\pi k}{2}\right) = (-1)^k$ when $k$ is replaced with $2k$, we obtain
\begin{equation}\label{Harm2k1}
\sum_{k=1}^{\infty} \frac{(-1)^{k} H_{2k}}{(2k+1)^2} = \Re\left(-\img\sum_{k=1}^\infty \frac{\eul^{\frac{\img \pi}{2}(k+1)} H_k}{(k+1)^2}\right).
\end{equation}
Substituting $z=\img$ in \eqref{lem6} and using \eqref{remark2}, we derive
\begin{equation}\label{Harm2k2}
\sum_{k=1}^{\infty} \frac{(-1)^{k} H_{2k}}{(2k+1)^2} = \Re(\img \li_3(1-\img)) - \frac{\G \ln{2}}{2} - \frac{\pi^3}{32} - \frac{\pi}{16}\ln^2{2}.
\end{equation}
Using Lemma \ref{lemconj1}, \eqref{Harm2k2} becomes
\begin{equation}\label{Harm2k3}
\sum_{k=1}^{\infty} \frac{(-1)^{k} H_{2k}}{(2k+1)^2} = \Im(\li_3(1+\img)) - \frac{\G \ln{2}}{2} - \frac{\pi^3}{32} - \frac{\pi}{16}\ln^2{2}.
\end{equation}
Utilizing the integral representation \eqref{harmonic1} for the harmonic numbers, we have
\begin{equation}\label{seria1}
\begin{split}
\sum_{k=1}^{\infty} \frac{(-1)^k H_k}{(2k+1)^2} &= \int_0^1 \int_0^u \int_0^1 \frac{-t^2}{u(1+t^2)(1+xt^2)}\, \rmd{d}x\, \rmd{d}t \, \rmd{d}u 
\\&= \int_0^1 \int_0^1 \frac{2\left(v \arctan{u} - \arctan(vu)\right)}{u(1- v^2)} \, \rmd{d}v\, \rmd{d}u
\\&= \int_0^1\int_0^1 \frac{2v\ln{u}}{1-v^2}\left(\frac{1}{1+ u^2 v^2} - \frac{1}{1+u^2}\right) \rmd{d}v \, \rmd{d}u = \int_0^1 \frac{\ln{u}\ln\left(1 + u^2\right)}{1+u^2}\, \rmd{d}u.
\end{split}
\end{equation}
By utilizing \eqref{lnlnsq134} and \eqref{Harm2k3}, we find
\begin{equation}\label{Harm2k4}
\sum_{k=1}^{\infty} \frac{(-1)^{k} H_{k}}{(2k+1)^2} = \int_0^1 \frac{\ln{u}\ln\left(1 + u^2\right)}{1+u^2}\, \rmd{d}u = -2\Im(\li_3(1+\img)) - \G \ln{2} + \frac{3\pi^3}{32} + \frac{\pi}{8}\ln^2{2}.
\end{equation}
The last series in \eqref{vuleanf} corresponds to the value of $\Theta_2(0, 1)$ in Theorem \ref{thm2}. Now, having found all the required values, \eqref{vuleanf} becomes
\begin{equation}\label{vuleanf2}
\begin{split}
\mathcal{V}_1 &= 4\G\ln{2} - \frac{\pi^3}{16}.
\end{split}
\end{equation}
\end{proof}

\begin{remark}
From Equations \eqref{Harm2k3} and \eqref{Harm2k4} we discover the following identity
\begin{equation}\label{h2khk}
\sum_{k=1}^{\infty} \frac{(-1)^k (H_k + 2H_{2k})}{(2k+1)^2} = -2\G \ln{2} + \frac{\pi^3}{32}.
\end{equation}
\textsf{Mathematica~13} does not offer a closed-form expression for \eqref{h2khk}, as we have presented, but it does provide a numerical approximation accurate to 15 decimal places.
\end{remark}

\subsection{Closed forms for a different digamma series}
In this subsection, we provide closed forms for a digamma series different from the class of series $\Theta_1(n, \alpha)$, $\Theta_2(n, \alpha)$ that are considering. We start with the logarithmic representation of $\arctan{t}$, which gives
\begin{equation}\label{anthmin1}
\int_0^1 \frac{\arctan{x} \ln(1- x)}{x} \, \rmd{d}x = -\frac{\img}{2} \int_0^1 \frac{\ln\left(1 + x^2\right) \ln(1- x)}{x} \, \rmd{d}x + \img \int_0^1 \frac{\ln(1- \img x) \ln(1 - x)}{x}\, \rmd{d}x.
\end{equation}
For the first integral, we already provided its closed form in \eqref{lnpnexsq}. Now, let's evaluate the second integral. Using the series expansion of $\ln(1-\img x)$ and then the beta function, the integral becomes
\begin{align*}
 \int_0^1 \frac{\ln(1- \img x) \ln(1 - x)}{x}\, \rmd{d}x &= -\sum_{k=1}^\infty \frac{\img^k}{k}\int_0^1 x^{k-1} \ln(1-x) \, \rmd{d}x = \sum_{k=1}^\infty \frac{\img^k H_k}{k^2}.
\end{align*}
Next, we rearrange \eqref{lem6} to obtain
\begin{equation}\label{lem6rear}
\sum_{k=1}^\infty \frac{H_k}{k^2} z^{k} = \zeta(3) + \li_2(1-z)\ln(1-z) + \li_3(z) - \li_3(1-z) + \frac{1}{2}\ln{z}\ln^2(1-z).
\end{equation}
Substituting $z=\img$ in \eqref{lem6rear}, and using \eqref{remark2} and \eqref{li3i}, we obtain
\begin{equation}\label{lemout1}
\sum_{k=1}^\infty \frac{H_k}{k^2} \img^{k} = \frac{29}{32}\zeta(3)  - \frac{\G\pi}{4} + \frac{\pi^2 \ln{2}}{32} - \img\left(\frac{\pi \ln^2{2}}{16}  + \frac{\G\ln{2}}{2}\right) - \li_3(1-\img).
\end{equation}
Applying Lemma \ref{lemconj1} to \eqref{lemli3ref}, we arrive at
\begin{equation}\label{li3pm1s}
\Re(\li_3(1\pm\img))  =-\frac{\pi^2}{32} \ln{2} + \frac{35}{64}\zeta(3).
\end{equation}
Using \eqref{li3pm1s} in \eqref{lemout1}, and noting that by Lemma \ref{lemconj1}, $\Im(\li_3(1\pm\img)) = -\Im(\li_3(1\mp\img))$, we simplify to
\begin{equation}\label{lemoutf}
\sum_{k=1}^\infty \frac{H_k}{k^2} \img^{k} = \frac{23}{64}\zeta(3) - \frac{\G\pi}{4} + \img\left(\Im(\li_3(1+\img)) - \frac{\G\ln{2}}{2} - \frac{\pi}{16}\ln^2{2}\right).
\end{equation}
Further using \eqref{lnpnexsq} in \eqref{anthmin1}, we get
\begin{equation}\label{anthmin2}
\int_0^1 \frac{\arctan{x} \ln(1- x)}{x} \, \rmd{d}x = -\Im(\li_3(1+\img)) + \frac{\G\ln{2}}{2} + \frac{\pi}{16}\ln^2{2}.
\end{equation}
Additionally, using the series expansion of $\ln(1-x)$ in \eqref{anthmin1}, we obtain
\begin{equation}\label{anthmin3}
\int_0^1 \frac{\arctan{x} \ln(1- x)}{x} \, \rmd{d}x = -1 + \frac{\pi^2}{6} - \frac{\pi^3}{24} - \frac{1}{4}\sum_{k=1}^\infty  \frac{\psi\left(\frac{k+5}{4}\right) - \psi\left(\frac{k+3}{4}\right)}{k^2}.
\end{equation}
By comparing \eqref{anthmin2} and \eqref{anthmin3}, we find
\begin{equation}\label{genr0508}
\sum_{k=1}^\infty  \frac{\psi\left(\frac{k+5}{4}\right) - \psi\left(\frac{k+3}{4}\right)}{k^2} = -4 + \frac{2\pi^2}{3} - \frac{\pi^3}{6} - 2\G\ln{2} - \frac{\pi}{4}\ln^2{2} + 4\Im(\li_3(1+\img)).
\end{equation}
In the following theorem, we provide a generalization of the series in \eqref{genr0508}.
\begin{theorem}Let $n \in \N$. Then \label{away1}
\begin{align*}
&\sum_{k=1}^\infty  \frac{\psi\left(\frac{k+2n+5}{4}\right) - \psi\left(\frac{k+2n+3}{4}\right)}{k^2}
\\&\quad= \begin{cases} 
\begin{aligned}[b]
-4 + \frac{2\pi^2}{3} - \frac{\pi^3}{6} - 2\G\ln{2} - \frac{\pi}{4}\ln^2{2} + 4\Im(\li_3(1+\img)),
\end{aligned} & \textup{if $n=0$,} \\ \ \\
\begin{aligned}[b]
\frac{86}{27} - \frac{4\pi^2}{9} + \frac{\pi^3}{6} + 2\G\ln{2} + \frac{\pi}{4}\ln^2{2} - 4\Im(\li_3(1+\img)),
\end{aligned} & \textup{if $n=1$},\\ \ \\
\begin{aligned}[b]&\frac{86}{27} - 2\sum_{j=1}^{\frac{n-1}{2}} \frac{8j+5}{(2j+1)(4j+3)^3}  + 32\sum_{j=1}^{\frac{n-1}{2}}  \frac{(2j+1)H_{4j+1}}{(4j+3)^2(4j+1)^2}
\\&- \left(\frac{4}{9} + \frac{4}{3} \sum_{j=1}^{\frac{n-1}{2}} \frac{1}{(4j+3)(4j+1)}\right)\pi^2 + \frac{\pi^3}{6} 
\\&+ 2\G\ln{2} + \frac{\pi}{4}\ln^2{2} - 4\Im(\li_3(1+\img)),
\end{aligned} & \begin{tabular}{cc}\textup{if $n$ is odd}\\\textup{and $n \neq 1$,}\end{tabular}\\ \ \\
\begin{aligned}[b]&-4 - \sum_{j=0}^{\frac{n-2}{2}} \frac{8j+9}{(j+1)(4j+5)^3}  + 64\sum_{j=0}^{\frac{n-2}{2}}  \frac{(j+1)H_{4j+3}}{(4j+5)^2(4j+3)^2}
\\&+ \left(\frac{2}{3} - \frac{4}{3} \sum_{j=0}^{\frac{n-2}{2}} \frac{1}{(4j+5)(4j+3)}\right)\pi^2 - \frac{\pi^3}{6} 
\\&-2\G\ln{2} - \frac{\pi}{4}\ln^2{2} + 4\Im(\li_3(1+\img)),
\end{aligned}& \textup{if $n$ is even.}
\end{cases}
\end{align*}
\end{theorem}
\begin{proof}
The proof of this theorem employs the same methodology as the proof of Theorem \ref{thm1}.
\end{proof}
\subsection{New Proofs for three more V\u{a}lean integrals} \label{sec3.4}
In this subsection, we present new proofs for three more integrals found in V\u{a}lean's book \cite[\S1.38, pp.~50, (1.174), (1.175), (1.177)]{bib9}. For ease of reference, we denote them as $\mathcal{V}_2$, $\mathcal{V}_3$, and $\mathcal{V}_4$, respectively. Now, let's begin with the first integral
\begin{equation}\label{val211}
\mathcal{V}_2 = \int_0^1 \frac{\arctan{x}\ln\left(1 + x^2\right)}{x} \, \rmd{d}x - 2\int_0^1 \frac{x\arctan{x}\ln\left(1 + x^2\right)}{1+x^2} \, \rmd{d}x + 2\int_0^1 \frac{x\arctan{x}\ln{x}}{1+x^2} \rmd{d}x.
\end{equation}
By applying integration by parts, we obtain
\begin{equation}\label{tobeca1}
2\int_0^1 \frac{x\arctan{x}\ln{x}}{1+x^2} \rmd{d}x = -\int_0^1 \frac{\ln\left(1+x^2\right) \ln{x}}{1+x^2} \, \rmd{d}x - \int_0^1 \frac{\ln\left(1+x^2\right) \arctan{x}}{x}\, \rmd{d}x.
\end{equation}
Substituting \eqref{tobeca1} into \eqref{val211}, we arrive at
\begin{equation}\label{val2f}
\mathcal{V}_2 = - 2\int_0^1 \frac{x\arctan{x}\ln\left(1 + x^2\right)}{1+x^2} \, \rmd{d}x  -\int_0^1 \frac{\ln\left(1+x^2\right) \ln{x}}{1+x^2} \, \rmd{d}x.
\end{equation}
For the first integral in \eqref{val2f}, we have by integration by parts
\begin{equation}\label{tontbeca1}
\int_0^1 \frac{2 x\arctan{x}\ln\left(1 + x^2\right)}{1+x^2} \, \rmd{d}x = \frac{\pi}{8}\ln^2{2} - \frac{1}{2}\int_0^1 \frac{\ln^2\left(1+x^2\right)}{1+x^2}\, \rmd{d}x.
\end{equation}
To evaluate the integral on the right side of \eqref{tontbeca1}, we utilize the trigonometric substitution $x=\tan\theta$ and use the result \eqref{cossq11s} to simplify \eqref{tontbeca1} as follows
\begin{equation}\label{vallast1}
\begin{split}
\int_0^1 \frac{2 x\arctan{x}\ln\left(1 + x^2\right)}{1+x^2} \, \rmd{d}x &= \frac{\pi}{8}\ln^2{2} - 2 \int_0^1 \ln^2(\cos{\theta}) \, \rmd{d}\theta
\\&= \frac{\pi}{8}\ln^2{2} + 2\G\ln{2} - \frac{\pi}{2} \ln^2{2} - \frac{\pi^3}{96} - 2\sum_{k=1}^{\infty} \frac{(-1)^{k-1} H_{2k}}{(2k+1)^2}
\\&= \G\ln{2} - \frac{7\pi^3}{96} - \frac{\pi}{2}\ln^2{2} +  2\Im(\li_3(1+\img)).
\end{split}
\end{equation}
Finally, using \eqref{lnlnsq134} in \eqref{val2f}, we find
\begin{align*}
\mathcal{V}_2 &= 2\G \ln{2} - \frac{\pi^3}{32} - 2 \sum_{k=1}^{\infty} \frac{(-1)^{k-1} H_{2k}}{(2k+1)^2} - \frac{\pi}{8}\ln^2{2} - 2\G\ln{2} + \frac{\pi}{2} \ln^2{2} 
\\&\quad+ \frac{\pi^3}{96} + 2\sum_{k=1}^{\infty} \frac{(-1)^{k-1} H_{2k}}{(2k+1)^2}
\\&= \frac{3\pi}{8}\ln^2{2} - \frac{\pi^3}{48}.
\end{align*}
Next, we substitute \eqref{arctansq1} and \eqref{Harm2k4} into \eqref{tobeca1} to get
\begin{equation}\label{atanlnexsq}
\begin{split}
2\int_0^1 \frac{x\arctan{x}\ln{x}}{1+x^2} \rmd{d}x &=  2\Im(\li_3(1+\img)) + \G \ln{2} - \frac{3\pi^3}{32} - \frac{\pi}{8}\ln^2{2}- \G \ln{2} - \frac{\pi^3}{16} 
\\&\quad - \frac{\pi}{8}\ln^2{2} + 2\Im(\li_3(1+\img)) 
\\&=  - \frac{5\pi^3}{32} - \frac{\pi}{4}\ln^2{2} + 4\Im(\li_3(1+\img)).
\end{split}
\end{equation}
From \eqref{Harm2k3} and \eqref{Harm2k4}, we derive
\begin{equation}\label{neweq1}
\sum_{k=1}^\infty \frac{(-1)^k\left(H_{2k} - H_k\right)}{(2k+1)^2} = 3\Im(\li_3(1+\img)) + \frac{\G \ln{2}}{2} - \frac{\pi^3}{8} - \frac{3\pi}{16}\ln^2{2}.
\end{equation}
Rearranging \eqref{neweq1}, we have
\begin{equation}\label{proofn1ab}
\sum_{k=1}^\infty \frac{(-1)^k\left(H_{2k+1} - H_k\right)}{(2k+1)^2} = 3\Im(\li_3(1+\img)) + \frac{\G \ln{2}}{2} - \frac{3 \pi^3}{32} - 1- \frac{3\pi}{16}\ln^2{2}.
\end{equation}
Using \eqref{proofn1ab} in \eqref{atanoneplx}, we obtain
\begin{equation}\label{proofn2ab}
\int_0^1 \frac{\arctan{x}\ln\left(1+ x\right)}{x} \rmd{d}x = \frac{3\G\ln{2}}{2} + \frac{3\pi^3}{32} + \frac{3\pi}{16}\ln^2{2} - 3\Im(\li_3(1+\img)).
\end{equation}
The third V\u{a}lean integral is given by
\begin{equation}\label{valeqn3}
\mathcal{V}_3 = \int_0^1 \frac{\arctan{x}\ln(1-x)}{x}\, \rmd{d}x + \int_0^1 \frac{\arctan{x}\ln(1+x)}{x} \, \rmd{d}x + 2\int_0^1 \frac{x\arctan{x}\ln{x}}{1+x^2} \rmd{d}x.
\end{equation}
Substituting \eqref{anthmin2}, \eqref{atanlnexsq} and \eqref{proofn2ab} into \eqref{valeqn3}, we obtain
\begin{align*}
 \mathcal{V}_3 &= 2\ln{2}\G - \frac{\pi^3}{16}.
\end{align*}
Moving on to the fourth V\u{a}lean integral
\begin{equation}\label{valeqn4}
\mathcal{V}_4 = \int_0^1 \frac{\arctan{x}\ln\left(1+x^2\right)}{x}\, \rmd{d}x + 2\int_0^1 \frac{x\arctan{x}\ln\left(1+x^2\right)}{1 + x^2} \,\rmd{d}x.
\end{equation}
Substituting \eqref{arctansq1} and \eqref{Harm2k3} into \eqref{valeqn4}, we obtain
\begin{equation}\label{valeanfour}
\begin{split}
\mathcal{V}_4  &= 2\G\ln{2} - \frac{3\pi}{8}\ln^2{2} - \frac{\pi^3}{96}.
\end{split}
\end{equation}
With these, we have now successfully proven all the integrals proposed by V\u{a}lean in \cite[\S1.38, pp.~50]{bib9}.

\subsection{New closed forms for some challenging integrals}\label{sec3.5}
In this subsection, building upon our established results, we present novel integrals that have not been previously explored in the literature. The first remarkable result is an integral representation of $\pi$, expressed as follows
\begin{equation}\label{newinegabd1}
\frac{32}{3}\int_0^1 \frac{\arctan{x} \ln(1+x)}{x} \, \rmd{d}x - 32 \int_0^1 \frac{\arctan{x} \ln(1-x)}{x} \, \rmd{d}x =  \pi^3.
\end{equation}

\begin{proof}
Equation \eqref{newinegabd1} is obtained by eliminating $\Im(\li_3(1+\img))$ from \eqref{anthmin2} and \eqref{proofn2ab}.
\end{proof}
Next, we present an integral that evaluates to $0$
\begin{equation}\label{newinega12} 
3\int_0^1 \frac{\arctan{x} \ln\left(1+x^2\right)}{x} \, \rmd{d}x - 2\int_0^1 \frac{\arctan{x} \ln\left(1+x\right)}{x} \, \rmd{d}x = 0.
\end{equation}
\begin{proof}
Equation \eqref{newinega12} is obtained by eliminating $\Im(\li_3(1+\img))$ from \eqref{arctansq1} and \eqref{proofn2ab}.
\end{proof}
Now, we introduce another integral that evaluates to $\pi$
\begin{equation}\label{newman2}
16\int_0^1 \frac{\arctan{x} \ln\left(1+x^2\right)}{x} \, \rmd{d}x - 32 \int_0^1 \frac{\arctan{x} \ln(1-x)}{x} \, \rmd{d}x =  \pi^3.
\end{equation}
\begin{proof}
Equation \eqref{newman2} is obtained by substituting \eqref{newinega12} into \eqref{newinegabd1}.
\end{proof}
The identity follows as the fourth integral
\begin{equation}\label{newinega123}
6\int_0^1 \frac{x \arctan{x} \ln{x}}{1+x^2} \, \rmd{d}x + 4\int_0^1 \frac{\arctan{x} \ln(1+x)}{x}\, \rmd{d}x = -\frac{3 \pi^3}{32} + 6\G \ln{2}.
\end{equation}
\begin{proof}
We derive \eqref{newinega123} by eliminating $\Im(\li_3(1+\img))$ from \eqref{atanlnexsq} and \eqref{proofn2ab}.
\end{proof}
Moving on, the fifth integral is
\begin{equation}\label{newman1}
\int_0^1 \frac{x \arctan{x} \ln{x}}{1+x^2} \, \rmd{d}x + \int_0^1 \frac{\arctan{x} \ln\left(1+x^2\right)}{x}\, \rmd{d}x = -\frac{\pi^3}{64} + \G \ln{2}.
\end{equation}
\begin{proof}
We obtain \eqref{newman1}  by substituting \eqref{newinega12} into \eqref{newinega123}.
\end{proof}
The sixth integral has the form
\begin{equation}\label{newman3}
\int_0^1 \frac{x \arctan{x} \ln{x}}{1+x^2} \, \rmd{d}x + 2\int_0^1 \frac{\arctan{x} \ln(1-x)}{x}\, \rmd{d}x = -\frac{5\pi^3}{64} + \G \ln{2}.
\end{equation}
\begin{proof}
We obtain \eqref{newman3}  by substituting \eqref{newman2} into \eqref{newman1}.
\end{proof}
The seventh integral is
\begin{equation}\label{newinega12345}
3\int_0^1 \frac{\ln{x}\ln\left(1+x^2\right)}{1+x^2}\, \rmd{d}x - 2\int_0^1 \frac{\arctan{x}\ln(1+x)}{x}\, \rmd{d}x = -6\G\ln{2} + \frac{3\pi^3}{32}.
\end{equation}
\begin{proof}
By eliminating $\Im(\li_3(1+\img))$ from the expressions in \eqref{Harm2k4} and \eqref{proofn2ab}, we establish \eqref{newinega12345}.
\end{proof}
The eighth integral is
\begin{equation}\label{newmanm1}
\int_0^1 \frac{\ln{x}\ln\left(1+x^2\right)}{1+x^2}\, \rmd{d}x - \int_0^1 \frac{\arctan{x}\ln\left(1+x^2\right)}{x}\, \rmd{d}x = -2\G\ln{2} + \frac{\pi^3}{32}.
\end{equation}
\begin{proof}
We obtain \eqref{newmanm1} by substituting \eqref{newinega12} into \eqref{newinega12345}.
\end{proof}
The ninth integral is
\begin{equation}\label{tenthint}
\int_0^1 \frac{\ln{x}\ln\left(1+x^2\right)}{1+x^2}\, \rmd{d}x - 2\int_0^1 \frac{\arctan{x}\ln(1-x)}{x}\, \rmd{d}x = -2\G\ln{2} + \frac{3\pi^3}{32}.
\end{equation}
\begin{proof}
We obtain \eqref{tenthint} by substituting \eqref{newinegabd1} into \eqref{newinega12345}.
\end{proof}
The tenth integral is
\begin{equation}\label{sixthint}
\int_0^1 \frac{\arctan{x}\ln(1+x)}{x}\, \rmd{d}x +  3\int_0^1 \frac{x\arctan{x}\ln\left(1+x^2\right)}{1+x^2} \, \rmd{d}x = 3\G\ln{2} - \frac{\pi^3}{64} - \frac{9\pi}{16}\ln^2{2}.
\end{equation}
\begin{proof}
We obtain \eqref{sixthint} by eliminating $\Im(\li_3(1+\img))$ from \eqref{vallast1} and \eqref{proofn2ab}.
\end{proof}
The eleventh integral is
\begin{equation}\label{t3nth}
\int_0^1 \frac{\arctan{x}\ln(1-x)}{x}\, \rmd{d}x + \int_0^1 \frac{x\arctan{x}\ln\left(1+x^2\right)}{1+x^2} \, \rmd{d}x = \G\ln{2} - \frac{7\pi^3}{192} - \frac{3\pi}{16}\ln^2{2}.
\end{equation}
\begin{proof}
We obtain \eqref{t3nth} by substituting \eqref{newinegabd1} into \eqref{sixthint}.
\end{proof}
The twelfth integral is
\begin{equation}\label{t4nth}
\int_0^1 \frac{\ln{x}\ln\left(1+x^2\right)}{1+x^2}\, \rmd{d}x + \int_0^1 \frac{x\arctan{x}\ln{x}}{1+x^2}\, \rmd{d}x = -\G\ln{2} + \frac{\pi^3}{64}.
\end{equation}
\begin{proof}
We obtain \eqref{t4nth} by adding \eqref{newman3} and \eqref{tenthint}.
\end{proof}
The thirteenth integral is
\begin{equation}\label{t5nth}
2\int_0^1 \frac{x\arctan{x}\ln\left(1+x^2\right)}{1+x^2} \, \rmd{d}x - \int_0^1 \frac{x\arctan{x}\ln{x}}{1+x^2}\, \rmd{d}x = \G\ln{2} + \frac{\pi^3}{192} - \frac{3\pi}{8}\ln^2{2}.
\end{equation}
\begin{proof}
Using \eqref{newman3} and \eqref{t3nth}, we obtain \eqref{t5nth}.
\end{proof}

\begin{remark}It is noteworthy that the integrals \eqref{newinegabd1}, \eqref{newinega12} and \eqref{newman2} were initially generalized by V\u{a}lean in  \cite[\S 1.37, (1.170)]{bib9}, \cite[Chapter 3, pp.~150-154]{bib9} and \cite[\S3.36, pp.~276]{bib2}, respectively. Nevertheless, our approach to these three integrals differs significantly from V\u{a}lean's, which adds value to our rediscovery. The remaining ten integrals, not found elsewhere in the literature, defy closed-form solutions even with \textsf{Mathematica~13}; it offers only numerical approximations. A remarkable observation is that by combining these integrals in various ways, one can derive even more new identities. For example, V\u{a}lean's fourth integral $\mathcal{V}_4$ can be obtained by substituting \eqref{newinega12} into \eqref{sixthint}.
\end{remark}

\subsection{Closed forms for $\Theta_1(n, 2m)$} \label{sec3.6}
In this subsection, we follow a similar approach to the proofs of Theorems \ref{thm1} and \ref{thm2} to derive the closed form of $\Theta_1(n, 2m)$. We begin by determining a closed form for the series $\Theta_1(0,0)$, which is given by
\begin{equation}\label{haf2}
\begin{split}
4\Theta_1(0,0) &= 4\sum_{k=1}^\infty \frac{1}{k^2} \sum_{j=0}^\infty \frac{(-1)^j}{2j+2k+3} = 4 \sum_{k=1}^\infty \frac{1}{k^2} \int_0^1 \frac{x^{2k+2}}{1+x^2} \, \rmd{d}x 
\\&= 4 \int_0^1 \li_2\left(x^2\right) - \frac{\li_2\left(x^2\right)}{1+x^2}\, \rmd{d}x  = \sum_{k=1}^{\infty} \frac{4}{k^2(2k+1)} -   \int_0^1 \frac{4\li_2\left(x^2\right)}{1+x^2}\, \rmd{d}x  
\\&= \frac{2\pi^2}{3} + 16\ln{2} - 16 - \frac{\pi^3}{6} - 8\int_0^1 \frac{\ln\left(1-x^2\right)}{x}\arctan{x}\, \rmd{d}x.
\end{split}
\end{equation}
By adding equations \eqref{anthmin2} and \eqref{proofn2ab}, we obtain
\begin{equation}\label{haf1}
\int_0^1 \frac{\ln\left(1-x^2\right)}{x}\arctan{x}\, \rmd{d}x = -4\Im(\li_3(1+\img)) + 2\G\ln{2} + \frac{3\pi^3}{32} + \frac{\pi}{4}\ln^2{2}.
\end{equation}
Substituting equation \eqref{haf1} into equation \eqref{haf2}, we find that
\begin{equation}\label{delta1zer}
\Theta_1(0,0) =  - 4\G\ln{2} + \frac{\pi^2}{6} + 4\ln{2} - 4 - \frac{11\pi^3}{48}  - \frac{\pi}{2} \ln^2{2} + 8\Im(\li_3(1+\img)).
\end{equation}
The theorem presented below provides a generalization from $\Theta_1(0,0)$ to $\Theta_1(n,0)$.
\begin{theorem}Let $n \in \N$. Then \label{thm6}
\begin{align*}
&\sum_{k=1}^\infty \frac{\psi\left(\frac{2k+2n+5}{4}\right) - \psi\left(\frac{2k+2n+3}{4}\right)}{(2k)^2}
\\&\quad= \begin{cases} 
\begin{aligned}[b]
- 4 + 4\ln{2} - \frac{\pi}{2} \ln^2{2} + \frac{\pi^2}{6} -  \frac{11\pi^3}{48} - 4\G\ln{2} + 8\Im(\li_3(1+\img)),
\end{aligned} & \textup{if $n=0$,} \\ \ \\
\begin{aligned}[b]
\frac{92}{27} - \frac{\pi^2}{9} - \frac{32\ln{2}}{9} + \frac{\pi}{2} \ln^2{2} + \frac{11\pi^3}{48} + 4\G\ln{2} - 8\Im(\li_3(1+\img)),
\end{aligned} & \textup{if $n=1$},\\ \ \\
\begin{aligned}[b]&\frac{92}{27} + 16\sum_{j=1}^{\frac{n-1}{2}} \frac{(2j+1)(2H_{4j+1} - H_{2j})}{(4j+3)^2(4j+1)^2} - 4\sum_{j=1}^{\frac{n-1}{2}} \frac{1}{(4j+3)^3}
\\&- \pi^2\left(\frac{1}{9} + \frac{1}{3}\sum_{j=1}^{\frac{n-1}{2}} \frac{1}{(4j+3)(4j+1)}\right)
\\&- \ln{2}\left(\frac{32}{9} + 32\sum_{j=1}^{\frac{n-1}{2}} \frac{2j+1}{(4j+3)^2(4j+1)^2}\right) + \frac{\pi}{2} \ln^2{2} 
\\&+\frac{11\pi^3}{48} + 4\G\ln{2} - 8\Im(\li_3(1+\img)), 
\end{aligned} & \begin{tabular}{cc}\textup{if $n$ is odd}\\ \ \textup{and $n \neq 1$,}\end{tabular}\\ \ \\
\begin{aligned}[b]&-4 + 32 \sum_{j=0}^{\frac{n-2}{2}} \frac{(j+1)\left(2H_{4j+3} - H_{2j+1}\right)}{(4j+5)^2(4j+3)^2} - 4\sum_{j=0}^{\frac{n-2}{2}} \frac{1}{(4j+5)^3}
\\&+ \pi^2\left(\frac{1}{6} - \frac{1}{3}\sum_{j=0}^{\frac{n-2}{2}} \frac{1}{(4j+5)(4j+3)}\right)
\end{aligned}
\end{cases}
\end{align*}
\newpage
\begin{align*}
\qquad\quad\begin{cases} 
\begin{aligned}[b]&+ \ln{2}\left(4 - 64\sum_{j=0}^{\frac{n-2}{2}} \frac{j+1}{(4j+5)^2(4j+3)^2}\right) - \frac{\pi}{2} \ln^2{2} 
\\&- \frac{11\pi^3}{48} - 4\G\ln{2} + 8\Im(\li_3(1+\img)), 
\end{aligned} & \textup{if $n$ is even.}
\end{cases}
\end{align*}
\end{theorem}
\begin{proof}
We use the same proof methodology as used in the proof of Theorem \ref{thm1}.
\end{proof}
We generalize $\Theta_1(n, 0)$ to $\Theta_1(n, 2m)$ for $m \in \N$ in the following theorem.
\begin{theorem}Let $m$ and $n$ be any two positive integers such that $n \geq m$. Then \label{thm7}
\begin{align*}
&\sum_{k=1}^\infty \frac{\psi\left(\frac{2k+2n+5}{4}\right) - \psi\left(\frac{2k+2n+3}{4}\right)}{(2k+2m)^2}
\\&\quad= \begin{cases} 
\begin{aligned}[b]
&- 4 - \sum_{k=1}^m \sum_{j=0}^k \frac{(-1)^{k+j}}{k^2(2j+1)} + 4\ln{2} - \left(\frac{1}{2} \ln^2{2} +\sum_{k=1}^m \frac{(-1)^{k-1}}{4k^2}\right)\pi 
\\&+ \frac{\pi^2}{6} -  \frac{11\pi^3}{48} - 4\G\ln{2} + 8\Im(\li_3(1+\img)),
\end{aligned} & \textup{if $n=m$,} \\ \ \\
\begin{aligned}[b]
&\frac{92}{27} + \sum_{k=1}^m \sum_{j=0}^{k+1} \frac{(-1)^{k+j}}{k^2(2j+1)}  + \left(\frac{1}{2} \ln^2{2} - \sum_{k=1}^m \frac{(-1)^{k}}{4k^2}\right)\pi 
\\&- \frac{\pi^2}{9} + \frac{11\pi^3}{48} + 4\G\ln{2}  - \frac{32\ln{2}}{9} - 8\Im(\li_3(1+\img)),
\end{aligned} & \textup{if $n=m+1$},\\ \ \\
\begin{aligned}[b]&\frac{92}{27} + 16\sum_{j=1}^{\frac{n-m-1}{2}} \frac{(2j+1)(2H_{4j+1} - H_{2j})}{(4j+3)^2(4j+1)^2} 
\\&- 4\sum_{j=1}^{\frac{n-m-1}{2}} \frac{1}{(4j+3)^3} -  \sum_{k=1}^m \sum_{j=0}^{k+n-m} \frac{(-1)^{j+k+n-m}}{k^2(2j+1)} 
\\&+ \left(\frac{1}{2} \ln^2{2}  + (-1)^{n-m} \sum_{k=1}^m \frac{(-1)^{k}}{4k^2}\right)\pi
\\&- \pi^2\left(\frac{1}{9} + \frac{1}{3}\sum_{j=1}^{\frac{n-m-1}{2}} \frac{1}{(4j+3)(4j+1)}\right)
\\&- \ln{2}\left(\frac{32}{9} + 32\sum_{j=1}^{\frac{n-m-1}{2}} \frac{2j+1}{(4j+3)^2(4j+1)^2}\right) 
\\&+\frac{11\pi^3}{48} + 4\G\ln{2} - 8\Im(\li_3(1+\img)), 
\end{aligned}&  \begin{tabular}{cc}\textup{if $n-m$ is odd}\\\textup{and $n-m \neq 1$.}\end{tabular}\\ \ \\
\begin{aligned}[b]
&-4 + 32 \sum_{j=0}^{\frac{n-m-2}{2}} \frac{(j+1)\left(2H_{4j+3} - H_{2j+1}\right)}{(4j+5)^2(4j+3)^2} 
\\&- 4\sum_{j=0}^{\frac{n-m-2}{2}} \frac{1}{(4j+5)^3} -  \sum_{k=1}^m \sum_{j=0}^{k+n-m} \frac{(-1)^{j+k+n-m}}{k^2(2j+1)} 
\end{aligned}
\end{cases}
\end{align*}
\newpage
\begin{align*}
\qquad \quad\begin{cases}
\begin{aligned}[b]&+ \left(-\frac{1}{2} \ln^2{2}  + (-1)^{n-m} \sum_{k=1}^m \frac{(-1)^{k}}{4k^2}\right)\pi
\\&+ \pi^2\left(\frac{1}{6} - \frac{1}{3}\sum_{j=0}^{\frac{n-m-2}{2}} \frac{1}{(4j+5)(4j+3)}\right)
\\&+ \ln{2}\left(4 - 64\sum_{j=0}^{\frac{n-m-2}{2}} \frac{j+1}{(4j+5)^2(4j+3)^2}\right) 
\\&- \frac{11\pi^3}{48} - 4\G\ln{2} + 8\Im(\li_3(1+\img)), 
\end{aligned} & \textup{if $n-m$ is even}.
\end{cases}
\end{align*}
\end{theorem}
\begin{proof}
This theorem is proved using the same method as employed in the proof of Theorem \ref{thm3}.
\end{proof}
\subsection{Closed forms for $\Theta_2(n, 2m)$} \label{sec3.7}
In this subsection, we utilize a similar approach to the proof of Theorems \ref{thm1}, \ref{thm2} to derive the closed form of $\Theta_2(n, 2m)$. We start by determining a closed form for the series $\Theta_2(0,0)$. We have
\begin{equation}\label{qot1}
\begin{split}
4\Theta_2(0,0) &= 4\sum_{k=1}^\infty \frac{(-1)^k}{k^2} \sum_{j=0}^\infty \frac{(-1)^j}{2j+2k+3} = 4 \sum_{k=1}^\infty \frac{(-1)^k}{k^2} \int_0^1 \frac{x^{2k+2}}{1+x^2} \, \rmd{d}x 
\\&= 4 \int_0^1 \li_2\left(-x^2\right) - \frac{\li_2\left(-x^2\right)}{1+x^2}\, \rmd{d}x  = \sum_{k=1}^{\infty} \frac{4(-1)^k}{k^2(2k+1)} -   \int_0^1 \frac{4\li_2\left(-x^2\right)}{1+x^2}\, \rmd{d}x  
\\&= -16 + 4\pi + 8\ln{2} - \frac{\pi^2}{3} + \frac{\pi^3}{12} - 8\int_0^1 \frac{\ln\left(1+x^2\right)}{x}\arctan{x}\, \rmd{d}x.
\end{split}
\end{equation}
Substituting \eqref{arctansq1} into \eqref{qot1}, we obtain
\begin{align*}
\Theta_2(0,0) &= -4 + \pi + 2\ln{2} - \frac{\pi}{4}\ln^2{2} - \frac{\pi^2}{12} - \frac{5\pi^3}{48} - 2\G \ln{2}  + 4\Im(\li_3(1+\img))
\end{align*}
In the following theorem, we establish the generalization of $\Theta_2(0, 0)$ to $\Theta_2(n, 0)$.
\begin{theorem}Let $n \in \N$. Then \label{thm8}
\begin{align*}
&\sum_{k=1}^\infty (-1)^k \frac{\psi\left(\frac{2k+2n+5}{4}\right) - \psi\left(\frac{2k+2n+3}{4}\right)}{(2k)^2}
\\&\quad= \begin{cases} 
\begin{aligned}[b]
-4 + \pi + 2\ln{2} - \frac{\pi}{4}\ln^2{2} - \frac{\pi^2}{12} - \frac{5\pi^3}{48} - 2\G \ln{2}  + 4\Im(\li_3(1+\img)),
\end{aligned} & \textup{if $n=0$,} \\ \ \\
\begin{aligned}[b]
&\frac{116}{27}- \frac{10\pi}{9} - \frac{16\ln{2}}{9} + \frac{\pi}{4}\ln^2{2} + \frac{\pi^2}{18} + \frac{5\pi^3}{48} + 2\G \ln{2}  
\\&- 4\Im(\li_3(1+\img)),
\end{aligned} & \textup{if $n=1$}, \\ \ \\
\begin{aligned}[b]
&\frac{116}{27} + 8\sum_{k=1}^{\frac{n-1}{2}} \sum_{j=0}^{2k} \frac{(-1)^j(16k^2 + 16k + 5)}{(2j+1)(4k+3)^2 (4k+1)^2} - 4 \sum_{j=1}^{\frac{n-1}{2}} \frac{1}{(4j+3)^3}
\\&+\left(-\frac{10}{9} - \sum_{j=1}^{\frac{n-1}{2}} \frac{2(16j^2 + 16j + 5)}{(4j+3)^2 (4j+1)^2} + \frac{1}{4}\ln^2{2}\right)\pi 
\end{aligned}
\end{cases}
\end{align*}

\begin{align*}
\qquad\quad\begin{cases}
\begin{aligned}[b]&-\left(\frac{1}{9} + \sum_{j=1}^{\frac{n-1}{2}} \frac{2j+1}{(4j+3)^2 (4j+1)^2}\right)16\ln{2}\\
&+ \left(\frac{1}{18} + \frac{1}{6}\sum_{j=1}^{\frac{n-1}{2}} \frac{1}{(4j+3)(4j+1)}\right)\pi^2 + \frac{5\pi^3}{48} 
\\&+ 2\G \ln{2}  - 4\Im(\li_3(1+\img)),
\end{aligned} &   \begin{tabular}{cc}\textup{if $n$ is odd}\\\textup{and $n \neq 1$,}\end{tabular}\\ \ \\
\begin{aligned}[b]&-4 - 8\sum_{k=0}^{\frac{n-2}{2}} \sum_{j=0}^{2k+1} \frac{(-1)^j(16k^2 + 32k + 17)}{(2j+1)(4k+5)^2 (4k+3)^2} 
\\&- 4 \sum_{j=0}^{\frac{n-2}{2}} \frac{1}{(4j+5)^3} +\left(1 + \sum_{j=0}^{\frac{n-2}{2}} \frac{2(16j^2 + 32j + 17)}{(4j+5)^2 (4j+3)^2} - \frac{1}{4}\ln^2{2}\right)\pi 
\\&+\left(1 - 16\sum_{j=0}^{\frac{n-2}{2}} \frac{j+1}{(4j+5)^2 (4j+3)^2}\right)2\ln{2}
\\&+ \left(-\frac{1}{2} + \sum_{j=0}^{\frac{n-2}{2}} \frac{1}{(4j+5)(4j+3)}\right)\frac{\pi^2}{6} - \frac{5\pi^3}{48} 
\\&-2\G \ln{2}  + 4\Im(\li_3(1+\img)),
\end{aligned} & \textup{if $n$ is even}.
\end{cases}
\end{align*}
\end{theorem}
\begin{proof}
The methodology utilized in the proof of Theorem \ref{thm2} is also applied in proving this theorem.
\end{proof}
In the following theorem, we extend the generalization from $\Theta_2(n, 0)$ to $\Theta_2(n, 2m)$ for $m \in \N$.
\begin{theorem}Let $m$ and $n$ be any two positive integers such that $n \geq m$. Then \label{thm9}
\begin{align*}
&\sum_{k=1}^\infty  (-1)^k \frac{\psi\left(\frac{2k+2n+5}{4}\right) - \psi\left(\frac{2k+2n+3}{4}\right)}{(2k+2m)^2}
\\&\quad= \begin{cases} 
\begin{aligned}[b]
&(-1)^m\left(-4 - \sum_{k=1}^m \sum_{j=0}^k \frac{(-1)^{j}}{k^2(2j+1)} + \left(1 + \sum_{k=1}^m \frac{1}{4k^2}\right)\pi + 2\ln{2} \right.
\\&\left.- \frac{\pi}{4}\ln^2{2} - \frac{\pi^2}{12} - \frac{5\pi^3}{48} - 2\G \ln{2}  + 4\Im(\li_3(1+\img))\right),
\end{aligned} & \textup{if $n=m$,}  \\ \ \\
\begin{aligned}[b]
&(-1)^m\left(\frac{116}{27} + \sum_{k=1}^m \sum_{j=0}^{k+1} \frac{(-1)^{j}}{k^2(2j+1)}  - \left(\frac{10}{9} + \sum_{k=1}^m \frac{1}{4k^2}\right)\pi \right.
\\&\left. -\frac{16\ln{2}}{9} + \frac{\pi}{4}\ln^2{2} + \frac{\pi^2}{18} + \frac{5\pi^3}{48} + 2\G \ln{2} - 4\Im(\li_3(1+\img))\right),
\end{aligned} & \textup{if $n=m+1$}, \\ \ \\
\begin{aligned}[b]
&(-1)^m\left(\frac{116}{27} - 4 \sum_{j=1}^{\frac{n-m-1}{2}} \frac{1}{(4j+3)^3} -  \sum_{k=1}^m \sum_{j=0}^{k+n-m} \frac{(-1)^{j+n-m}}{k^2(2j+1)} \right.
\\&+ 8\sum_{k=1}^{\frac{n-m-1}{2}} \sum_{j=0}^{2k} \frac{(-1)^j(16k^2 + 16k + 5)}{(2j+1)(4k+3)^2 (4k+1)^2} 
\end{aligned} 
\end{cases}
\end{align*}
\newpage
\begin{align*}
\qquad\begin{cases}
\begin{aligned}[b]
&+ \left(-\frac{10}{9} - \sum_{j=1}^{\frac{n-m-1}{2}} \frac{2(16j^2 + 16j + 5)}{(4j+3)^2 (4j+1)^2} + \frac{1}{4}\ln^2{2} +  \sum_{k=1}^m \frac{(-1)^{n-m}}{4k^2}\right)\pi
\\&-\left(\frac{1}{9} + \sum_{j=1}^{\frac{n-m-1}{2}} \frac{2j+1}{(4j+3)^2 (4j+1)^2}\right)16\ln{2}
\\&+ \left(\frac{1}{18} + \frac{1}{6}\sum_{j=1}^{\frac{n-m-1}{2}} \frac{1}{(4j+3)(4j+1)}\right)\pi^2 + \frac{5\pi^3}{48} 
\\&\left.+ 2\G \ln{2}  - 4\Im(\li_3(1+\img))\right),
\end{aligned} &   \begin{tabular}{cc}\textup{if $n-m$ is odd}\\\textup{and $n-m \neq 1$,}\end{tabular}\\ \ \\
\begin{aligned}[b]&(-1)^m\left(-4  - 4 \sum_{j=0}^{\frac{n-m-2}{2}} \frac{1}{(4j+5)^3} -  \sum_{k=1}^m \sum_{j=0}^{k+n-m} \frac{(-1)^{j+n-m}}{k^2(2j+1)}\right.
\\&- 8\sum_{k=0}^{\frac{n-m-2}{2}} \sum_{j=0}^{2k+1} \frac{(-1)^j(16k^2 + 32k + 17)}{(2j+1)(4k+5)^2 (4k+3)^2} 
\\&+\left(1 + \sum_{j=0}^{\frac{n-m-2}{2}} \frac{2(16j^2 + 32j + 17)}{(4j+5)^2 (4j+3)^2} - \frac{1}{4}\ln^2{2} +  \sum_{k=1}^m \frac{(-1)^{n-m}}{4k^2}\right)\pi 
\\&+\left(1 - 16\sum_{j=0}^{\frac{n-m-2}{2}} \frac{j+1}{(4j+5)^2 (4j+3)^2}\right)2\ln{2}
\\&+ \left(-\frac{1}{2} + \sum_{j=0}^{\frac{n-m-2}{2}} \frac{1}{(4j+5)(4j+3)}\right)\frac{\pi^2}{6} - \frac{5\pi^3}{48} 
\\&\left.-2\G \ln{2}  + 4\Im(\li_3(1+\img))\right),
\end{aligned} & \textup{if $n-m$ is even.}
\end{cases}
\end{align*}
\end{theorem}
\begin{proof}
This theorem is proved using the same method as employed in the proof of Theorem \ref{thm4}.
\end{proof}

\subsection{New identities for some generalized digamma series}\label{sec3.8}
In the following theorems, we present several infinite series that result in closed forms of the type $a_0 + a_1 \pi + a_2 \pi^2 + a_3 \pi^3$, where the constants $a_0$, $a_1$, $a_2$, and $a_3$ are real-valued and determined by our formulas.
\begin{theorem}Let $n \in \N$. Then \label{thm10}
\begin{align*}
&\sum_{k=1}^\infty \frac{\left(\psi\left(\frac{2k+2n+5}{4}\right) - \psi\left(\frac{2k+2n+3}{4}\right)\right)\left(\frac12 - (-1)^k\right)}{(2k)^2}
\\&\quad= \begin{cases} 
\begin{aligned}[b]
2 - \pi + \frac{\pi^2}{6} - \frac{\pi^3}{96},
\end{aligned} & \textup{if $n=0$,} \\ \ \\
\begin{aligned}[b]
-\frac{70}{27} + \frac{10\pi}{9} - \frac{\pi^2}{9} + \frac{\pi^3}{96},
\end{aligned} & \textup{if $n=1$}, \\ \ \\
\begin{aligned}[b]
-\frac{70}{27} + 8\sum_{j=1}^{\frac{n-1}{2}} \frac{(2j+1)(2H_{4j+1} - H_{2j})}{(4j+3)^2(4j+1)^2} + 2\sum_{j=1}^{\frac{n-1}{2}} \frac{1}{(4j+3)^3}
\end{aligned} 
\end{cases}
\end{align*}

\newpage
\begin{align*}
 \qquad\begin{cases} 
\begin{aligned}[b]&- 8\sum_{k=1}^{\frac{n-1}{2}} \sum_{j=0}^{2k} \frac{(-1)^j(16k^2 + 16k + 5)}{(2j+1)(4k+3)^2 (4k+1)^2} 
\\&+ \left(\frac{10}{9} + \sum_{j=1}^{\frac{n-1}{2}} \frac{2(16j^2 + 16j + 5)}{(4j+3)^2(4j+1)^2}\right)\pi 
\\&- \left(\frac{1}{3} + \sum_{j=1}^{\frac{n-1}{2}} \frac{1}{(4j+3)(4j+1)}\right)\frac{\pi^2}{3} + \frac{\pi^3}{96},
\end{aligned} & \textup{if $n$ is odd and $n \neq 1$}, \\ \ \\
\begin{aligned}[b]&2 + 16 \sum_{j=0}^{\frac{n-2}{2}} \frac{(j+1)(2H_{4j+3} - H_{2j+1})}{(4j+5)^2(4j+3)^2} + 2\sum_{j=0}^{\frac{n-2}{2}} \frac{1}{(4j+5)^3}
\\&+ 8\sum_{k=0}^{\frac{n-2}{2}} \sum_{j=0}^{2k+1} \frac{(-1)^j(16k^2 + 32k + 17)}{(2j+1)(4k+5)^2 (4k+3)^2} 
\\& -\left(1 + \sum_{j=0}^{\frac{n-2}{2}} \frac{2(16j^2 + 32j + 17)}{(4j+5)^2(4j+3)^2}\right)\pi 
\\&+ \left(\frac{1}{2} - \sum_{j=0}^{\frac{n-2}{2}} \frac{1}{(4j+5)(4j+3)}\right) \frac{\pi^2}{3} - \frac{\pi^3}{96},
\end{aligned} & \textup{if $n$ is even.}
\end{cases}
\end{align*}
\end{theorem}
\begin{proof}
The proof is completed by subtracting by subtracting the entries of Theorem \ref{thm8} from half the entries of Theorem \ref{thm6}.
\end{proof}

\begin{ex} For $n=2$, we have
\begin{equation}
\sum_{k=1}^\infty \frac{\left(\psi\left(\frac{2k+9}{4}\right) - \psi\left(\frac{2k+7}{4}\right)\right)\left(\frac12 - (-1)^k\right)}{(2k)^2} = \frac{8804}{3375}-\frac{259 \pi }{225}+\frac{13 \pi^2}{90}-\frac{\pi^3}{96}.
\end{equation}
For $n=3$, we have
\begin{equation}
\sum_{k=1}^\infty \frac{\left(\psi\left(\frac{2k+11}{4}\right) - \psi\left(\frac{2k+9}{4}\right)\right)\left(\frac12 - (-1)^k\right)}{(2k)^2} =-\frac{3167372}{1157625}+\frac{12916 \pi }{11025}-\frac{38 \pi^2}{315}+\frac{\pi^3}{96}.
\end{equation}

For $n=4$, we have
\begin{equation}
\sum_{k=1}^\infty \frac{\left(\psi\left(\frac{2k+13}{4}\right) - \psi\left(\frac{2k+11}{4}\right)\right)\left(\frac12 - (-1)^k\right)}{(2k)^2} = \frac{85428394}{31255875}-\frac{117469 \pi }{99225}+\frac{263 \pi^2}{1890}-\frac{\pi^3}{96}.
\end{equation}
\end{ex}\noindent
By generalizing Theorem \ref{thm10} to two variables, we arrive at the following theorem.
\begin{theorem}Let $m$ and $n$ be any two positive integers such that $n \geq m$. Then \label{thm11}
\begin{align*}
&\sum_{k=1}^\infty  \frac{\left(\psi\left(\frac{2k+2n+5}{4}\right) - \psi\left(\frac{2k+2n+3}{4}\right)\right)\left(\frac{1}{2} - (-1)^k\right)}{(2k+4m)^2}
\\&\quad= \begin{cases} 
\begin{aligned}[b]
&2 - \sum_{k=1}^{2m} \sum_{j=0}^k \frac{(-1)^{j}\left(\frac{1}{2}(-1)^k - 1\right)}{k^2(2j+1)} + \left(-1 + \sum_{k=1}^{2m} \frac{\left(\frac{1}{2}(-1)^k - 1\right)}{4k^2}\right)\pi
\\&+ \frac{\pi^2}{6} - \frac{\pi^3}{96},
\end{aligned} & \textup{if $n=2m$,} 
\end{cases}
\end{align*}

\newpage
\begin{align*}
\begin{cases} 
\begin{aligned}[b]
&\frac{-70}{27} + \sum_{k=1}^{2m} \sum_{j=0}^{k+1} \frac{(-1)^{j}\left(\frac{1}{2}(-1)^k-1\right)}{k^2(2j+1)}  + \left(\frac{10}{9} - \sum_{k=1}^{2m} \frac{\left(\frac{1}{2}(-1)^k-1\right)}{4k^2}\right)\pi 
\\&-\frac{\pi^2}{9} + \frac{\pi^3}{96},
\end{aligned} & \textup{if $n=2m+1$},\\ \ \\
\begin{aligned}[b]&\frac{-70}{27} + 8\sum_{j=1}^{\frac{n-2m-1}{2}} \frac{(2j+1)(2H_{4j+1} - H_{2j})}{(4j+3)^2(4j+1)^2} 
\\&+ 2\sum_{j=1}^{\frac{n-2m-1}{2}} \frac{1}{(4j+3)^3}
 -  \sum_{k=1}^{2m} \sum_{j=0}^{k+n-2m} \frac{(-1)^{j+n}\left(\frac{1}{2}(-1)^k-1\right)}{k^2(2j+1)}
\\& -8\sum_{k=1}^{\frac{n-2m-1}{2}} \sum_{j=0}^{2k} \frac{(-1)^j(16k^2 + 16k + 5)}{(2j+1)(4k+3)^2 (4k+1)^2} 
\\&+ \left(\frac{10}{9} + \sum_{j=1}^{\frac{n-2m-1}{2}} \frac{2(16j^2 + 16j + 5)}{(4j+3)^2 (4j+1)^2} +  \sum_{k=1}^{2m} \frac{(-1)^n \left(\frac{1}{2}(-1)^k-1\right)}{4k^2}\right)\pi
\\&- \left(\frac{1}{9} + \frac{1}{3}\sum_{j=1}^{\frac{n-2m-1}{2}} \frac{1}{(4j+3)(4j+1)}\right)\pi^2 + \frac{\pi^3}{96} ,
\end{aligned} &   \begin{tabular}{cc}\textup{if $n-2m$ is odd}\\\textup{and $n-2m \neq 1$,}\end{tabular}\\ \ \\
\begin{aligned}[b]&2 + 16\sum_{j=1}^{\frac{n-2m-2}{2}} \frac{(j+1)(2H_{4j+3} - H_{2j+1})}{(4j+5)^2(4j+3)^2} 
\\&+2 \sum_{j=0}^{\frac{n-2m-2}{2}} \frac{1}{(4j+5)^3}  -  \sum_{k=1}^{2m} \sum_{j=0}^{k+n-2m} \frac{(-1)^{j+n}\left(\frac{1}{2}(-1)^k-1\right)}{k^2(2j+1)}
\\&+8\sum_{k=0}^{\frac{n-2m-2}{2}} \sum_{j=0}^{2k+1} \frac{(-1)^j(16k^2 + 32k + 17)}{(2j+1)(4k+5)^2 (4k+3)^2} 
\\&+\left(-1 - \sum_{j=0}^{\frac{n-2m-2}{2}} \frac{2(16j^2 + 32j + 17)}{(4j+5)^2 (4j+3)^2} +  \sum_{k=1}^{2m} \frac{(-1)^n\left(\frac{1}{2}(-1)^k-1\right)}{4k^2}\right)\pi 
\\&+ \left(\frac{1}{2} - \sum_{j=0}^{\frac{n-2m-2}{2}} \frac{1}{(4j+5)(4j+3)}\right)\frac{\pi^2}{3} - \frac{\pi^3}{96},
\end{aligned} & \textup{if $n-2m$ is even.}
\end{cases}
\end{align*}
\end{theorem}
\begin{proof}
By substituting $n=2m$ in the entries of Theorem \ref{thm7} and \ref{thm9}, and then subtracting the latter result of Theorem \ref{thm9} from half the latter result of Theorem \ref{thm7}, the proof is complete.
\end{proof}

\begin{ex}
For $n=6$ and $m=3$, we have
\begin{equation}
\sum_{k=1}^\infty  \frac{\left(\psi\left(\frac{2k+17}{4}\right) - \psi\left(\frac{2k+15}{4}\right)\right)\left(\frac{1}{2} - (-1)^k\right)}{(2k+12)^2} =\frac{1073869873}{324324000}-\frac{42457 \pi }{28800}+\frac{\pi ^2}{6}-\frac{\pi ^3}{96}.
\end{equation}
For $n=7$ and $m=3$, we have
\begin{equation}
\sum_{k=1}^\infty  \frac{\left(\psi\left(\frac{2k+19}{4}\right) - \psi\left(\frac{2k+17}{4}\right)\right)\left(\frac{1}{2} - (-1)^k\right)}{(2k+12)^2} = -\frac{681924389}{162162000}+\frac{5073\pi}{3200}-\frac{\pi^2}{9}+\frac{\pi^3}{96}.
\end{equation}
For $n=6$ and $m=2$, we have
\begin{equation}
\sum_{k=1}^\infty  \frac{\left(\psi\left(\frac{2k+17}{4}\right) - \psi\left(\frac{2k+15}{4}\right)\right)\left(\frac{1}{2} - (-1)^k\right)}{(2k+8)^2} = \frac{6775331}{1716000}-\frac{46277 \pi }{28800}+\frac{13 \pi^2}{90}-\frac{\pi^3}{96}.
\end{equation}
For $n=7$ and $m=2$, we have
\begin{equation}
\sum_{k=1}^\infty  \frac{\left(\psi\left(\frac{2k+19}{4}\right) - \psi\left(\frac{2k+17}{4}\right)\right)\left(\frac{1}{2} - (-1)^k\right)}{(2k+8)^2} = -\frac{78022319}{18393375}+\frac{2296373 \pi }{1411200}-\frac{38 \pi^2}{315}+\frac{\pi^3}{96}.
\end{equation}
\end{ex}
In the next theorem, we introduce another generalized digamma series that yields the closed form $a_0 + a_1 \pi + a_2 \pi^2 + a_3 \pi^3$, with $a_0,\ldots,a_3 \in \mathbb{Q}$. In these series, $1/2 - (-1)^k$ in Theorem \ref{thm11} is replaced with $1/2 + (-1)^k$.
\begin{theorem}Let $m$ and $n$ be any two positive integers such that $n \geq m$. Then \label{thm12}
\begin{align*}
&\sum_{k=1}^\infty  \frac{\left(\psi\left(\frac{2k+2n+5}{4}\right) - \psi\left(\frac{2k+2n+3}{4}\right)\right)\left(\frac{1}{2} + (-1)^k\right)}{(2k+4m-2)^2}
\\&\hspace{0.15cm}= \begin{cases} 
\begin{aligned}[b]
&2 - \sum_{k=1}^{2m-1} \sum_{j=0}^k \frac{(-1)^{j}\left(\frac{1}{2}(-1)^k - 1\right)}{k^2(2j+1)} 
\\&+ \left(-1 + \sum_{k=1}^{2m-1} \frac{\left(\frac{1}{2}(-1)^k - 1\right)}{4k^2}\right)\pi + \frac{\pi^2}{6} - \frac{\pi^3}{96},
\end{aligned} & \textup{if $n=2m-1$,} \\ \ \\
\begin{aligned}[b]
&-\frac{70}{27} + \sum_{k=1}^{2m-1} \sum_{j=0}^{k+1} \frac{(-1)^{j}\left(\frac{1}{2}(-1)^k-1\right)}{k^2(2j+1)}  
\\&+ \left(\frac{10}{9} - \sum_{k=1}^{2m-1} \frac{\left(\frac{1}{2}(-1)^k-1\right)}{4k^2}\right)\pi -\frac{\pi^2}{9} + \frac{\pi^3}{96},
\end{aligned} & \textup{if $n=2m$},\\ \ \\
\begin{aligned}[b]&-\frac{70}{27} + 8\sum_{j=1}^{\frac{n-2m}{2}} \frac{(2j+1)(2H_{4j+1} - H_{2j})}{(4j+3)^2(4j+1)^2} 
\\&+ 2\sum_{j=1}^{\frac{n-2m}{2}} \frac{1}{(4j+3)^3}
 +\sum_{k=1}^{2m-1} \sum_{j=0}^{k+n-2m+1} \frac{(-1)^{j+n}\left(\frac{1}{2}(-1)^k-1\right)}{k^2(2j+1)}
\\& -8\sum_{k=1}^{\frac{n-2m}{2}} \sum_{j=0}^{2k} \frac{(-1)^j(16k^2 + 16k + 5)}{(2j+1)(4k+3)^2 (4k+1)^2} 
\\&+ \left(\frac{10}{9} + \sum_{j=1}^{\frac{n-2m}{2}} \frac{2(16j^2 + 16j + 5)}{(4j+3)^2 (4j+1)^2} - \sum_{k=1}^{2m-1} \frac{(-1)^n \left(\frac{1}{2}(-1)^k - 1\right)}{4k^2}\right)\pi
\\&- \left(\frac{1}{9} + \frac{1}{3}\sum_{j=1}^{\frac{n-2m}{2}} \frac{1}{(4j+3)(4j+1)}\right)\pi^2 + \frac{\pi^3}{96} ,
\end{aligned} &   \begin{tabular}{cc}\textup{if}\\ \textup{$n-2m+1$}\\ \textup{is odd and}\\\ \textup{$n-2m+1 \neq 1$,}\end{tabular} \\ \ \\
\begin{aligned}[b]&2 + 16\sum_{j=1}^{\frac{n-2m-1}{2}} \frac{(j+1)(2H_{4j+3} - H_{2j+1})}{(4j+5)^2(4j+3)^2} 
\end{aligned}
\end{cases}
\end{align*}
\newpage
\begin{align*}
\quad\begin{cases}
\begin{aligned}[b]
&+2 \sum_{j=0}^{\frac{n-2m-1}{2}} \frac{1}{(4j+5)^3}  +  \sum_{k=1}^{2m-1} \sum_{j=0}^{k+n-2m+1} \frac{(-1)^{j+n}\left(\frac{1}{2}(-1)^k-1\right)}{k^2(2j+1)}\\
&+8\sum_{k=0}^{\frac{n-2m-1}{2}} \sum_{j=0}^{2k+1} \frac{(-1)^j(16k^2 + 32k + 17)}{(2j+1)(4k+5)^2 (4k+3)^2} 
\\&+\left(-1 - \sum_{j=0}^{\frac{n-2m-1}{2}} \frac{2(16j^2 + 32j + 17)}{(4j+5)^2 (4j+3)^2} -  \sum_{k=1}^{2m-1} \frac{(-1)^n\left(\frac{1}{2}(-1)^k-1\right)}{4k^2}\right)\pi 
\\&+ \left(\frac{1}{2} - \sum_{j=0}^{\frac{n-2m-1}{2}} \frac{1}{(4j+5)(4j+3)}\right)\frac{\pi^2}{3} - \frac{\pi^3}{96},
\end{aligned} & \begin{tabular}{cc}\textup{if} \\\textup{$n-2m+1$}\\\textup{is even}.\end{tabular}
\end{cases}
\end{align*}
\end{theorem}
\begin{proof}
By substituting $n=2m-1$ in the entries of Theorem \ref{thm7} and \ref{thm9}, and then adding the latter result of Theorem \ref{thm9} with half the latter result of Theorem \ref{thm7}, the proof is complete.
\end{proof}

\begin{ex}For $n=1$ and $m=1$, we have
\begin{equation}
\sum_{k=1}^\infty  \frac{\left(\psi\left(\frac{2k+7}{4}\right) - \psi\left(\frac{2k+5}{4}\right)\right)\left(\frac{1}{2} + (-1)^k\right)}{(2k+2)^2} = 3-\frac{11 \pi }{8}+\frac{\pi^2}{6}-\frac{\pi^3}{96}.
\end{equation}
For $n=2$ and $m=1$, we have
\begin{equation}
\sum_{k=1}^\infty  \frac{\left(\psi\left(\frac{2k+9}{4}\right) - \psi\left(\frac{2k+7}{4}\right)\right)\left(\frac{1}{2} + (-1)^k\right)}{(2k+2)^2} = -\frac{1051}{270}+\frac{107 \pi }{72}-\frac{\pi ^2}{9}+\frac{\pi ^3}{96}.
\end{equation}
For $n=4$ and $m=1$, we have
\begin{equation}
\sum_{k=1}^\infty  \frac{\left(\psi\left(\frac{2k+13}{4}\right) - \psi\left(\frac{2k+11}{4}\right)\right)\left(\frac{1}{2} + (-1)^k\right)}{(2k+2)^2} = -\frac{9234319}{2315250}+\frac{136403 \pi }{88200}-\frac{38 \pi^2}{315}+\frac{\pi^3}{96}.
\end{equation}
For $n=5$ and $m=1$, we have
\begin{equation}
\sum_{k=1}^\infty  \frac{\left(\psi\left(\frac{2k+15}{4}\right) - \psi\left(\frac{2k+13}{4}\right)\right)\left(\frac{1}{2} + (-1)^k\right)}{(2k+2)^2} = \frac{1323415409}{343814625}-\frac{1237427 \pi }{793800}+\frac{263 \pi^2}{1890}-\frac{\pi^3}{96}.
\end{equation}
\end{ex}
\begin{remark}
\textsf{Mathematica~13} does not produce closed forms for any value of $n$ in Theorem \ref{thm10}, nor for combinations of $n$ and $m$ values in Theorems \ref{thm11} and \ref{thm12}.  However, it does provide numerical approximations. Notably, in the three theorems the numerical approximation is accurate to around 12 decimal places.
\end{remark}
In the following theorem, we present two generalized digamma series also characterized by the closed form $a_0 + a_1 \pi + a_2 \pi^2 + a_3 \pi^3$. However, the coefficient $a_0$ in this series is a transcendental number, and $a_1=0$. More precisely, $a_0$ is of the form $b_0 + b_2 \ln{2}$, where $b_0$ and $b_2$ are rational numbers.
\begin{theorem}Let $n \in \N$. Then \label{thm13}
\begin{align*}
&\sum_{k=1}^\infty \frac{\frac{1}{8}\psi\left(\frac{2k+2n+5}{4}\right) - \psi\left(\frac{k+2n+5}{4}\right) + \psi\left(\frac{k+2n+3}{4}\right) - \frac{1}{8}\psi\left(\frac{2k+2n+3}{4}\right)}{k^2}
\\&\quad= \begin{cases} 
\begin{aligned}[b]
2(1 + \ln{2}) - \frac{7\pi^2}{12} + \frac{5\pi^3}{96},
\end{aligned} & \textup{if $n=0$,} \\ \ \\
\begin{aligned}[b]
-\frac{8}{9}\left(\frac{5}{3} + 2\ln{2}\right) + \frac{7\pi^2}{18} - \frac{5\pi^3}{96},
\end{aligned} & \textup{if $n=1$},\\ \ \\
\begin{aligned}[b]&\frac{-40}{27} - 8\sum_{j=1}^{\frac{n-1}{2}} \frac{(2j+1)(2H_{4j+1} + H_{2j})}{(4j+3)^2(4j+1)^2} + 4\sum_{j=1}^{\frac{n-1}{2}} \frac{3j+2}{(2j+1)(4j+3)^3}
\\&+ \frac{7\pi^2}{2}\left(\frac{1}{9} + \frac{1}{3}\sum_{j=1}^{\frac{n-1}{2}} \frac{1}{(4j+3)(4j+1)}\right)
\\&- 16\ln{2}\left(\frac{1}{9} + \sum_{j=1}^{\frac{n-1}{2}} \frac{2j+1}{(4j+3)^2(4j+1)^2}\right) -\frac{5\pi^3}{96},
\end{aligned} & \begin{tabular}{cc}\textup{if $n$ is odd}\\ \ \textup{and $n \neq 1$,}\end{tabular}\\ \ \\
\begin{aligned}[b]&2 - 16 \sum_{j=0}^{\frac{n-2}{2}} \frac{(j+1)\left(2H_{4j+3} + H_{2j+1}\right)}{(4j+5)^2(4j+3)^2} + \sum_{j=0}^{\frac{n-2}{2}} \frac{6j+7}{(j+1)(4j+5)^3}
\\&- \frac{7\pi^2}{2}\left(\frac{1}{6} - \frac{1}{3}\sum_{j=0}^{\frac{n-2}{2}} \frac{1}{(4j+5)(4j+3)}\right)
\\&+ 2\ln{2}\left(1 - 16\sum_{j=0}^{\frac{n-2}{2}} \frac{j+1}{(4j+5)^2(4j+3)^2}\right) + \frac{5\pi^3}{96},
\end{aligned}& \textup{if $n$ is even.}
\end{cases}
\end{align*}
\end{theorem}
\begin{proof}
By subtracting the results of Theorem \ref{away1} from half the results of Theorem \ref{thm6}, the proof is complete.
\end{proof}
\begin{ex}
For $n=2$, we have
\begin{equation}
\begin{split}
&\sum_{k=1}^\infty \frac{\frac{1}{8}\psi\left(\frac{2k+9}{4}\right) - \psi\left(\frac{k+9}{4}\right) + \psi\left(\frac{k+7}{4}\right) - \frac{1}{8}\psi\left(\frac{2k+7}{4}\right)}{k^2} 
\\&\qquad\qquad\qquad= \frac{11(529+570 \ln{2})}{3375} -\frac{91\pi^2}{180}  + \frac{5 \pi ^3}{96}.
\end{split}
\end{equation}
For $n=3$, we have
\begin{equation}
\begin{split}
&\sum_{k=1}^\infty \frac{\frac{1}{8}\psi\left(\frac{2k+11}{4}\right) - \psi\left(\frac{k+11}{4}\right) + \psi\left(\frac{k+9}{4}\right) - \frac{1}{8}\psi\left(\frac{2k+9}{4}\right)}{k^2} 
\\&\qquad\qquad\qquad= -\frac{4(457523+525840 \ln{2}) }{1157625}+ \frac{19 \pi^2}{45}-\frac{5 \pi^3}{96}
\end{split}
\end{equation}
\end{ex}

\begin{remark}
\textsf{Mathematica~13} does not yield closed-form solutions for any value of $n$ in Theorem \ref{thm13} either. It is also worth noting that the provided numerical approximation maintains precision up to 26 decimal places.
\end{remark}

\section{Conclusion}
In this study, we have demonstrated the derivation of closed forms for a variety of generalized digamma series and established their intriguing connection with the elegant expression $a_0 + a_1 \pi + a_2 \pi^2 + a_3 \pi^3$. Furthermore, our investigation has led us to unveil new closed forms for several integrals. Remarkably, we have also rediscovered two integral representations of $\pi$, which take on the exquisite forms
\begin{equation*}
\begin{split}
\pi &= \left(\int_0^1 \frac{\arctan{x}}{x} \left(2^4\ln\left(1+x^2\right)  - 2^5 \ln(1-x)\right) \rmd{d}x\right)^{\frac{1}{3}}
\\&= \left(\int_0^1 \frac{2^5 \arctan{x}}{x} \left(\frac{1}{3}\ln(1+x) -\ln(1-x)\right)  \rmd{d}x\right)^{\frac{1}{3}}.
\end{split}
\end{equation*}
Furthermore, from equations \eqref{Harm2k3} and \eqref{Harm2k4}, one can derive alternative infinite series representations for $\Im(\li_3(1+\img))$
\begin{align*}
\Im(\li_3(1+\img)) &= \frac{\G\ln{2}}{2} + \frac{\pi}{16}\ln^2{2} + \sum_{k=0}^\infty \frac{H_{2k+1}(-1)^k}{(2k+1)^2}
\\&= -\frac{\G\ln{2}}{2} + \frac{3\pi^3}{64} + \frac{\pi}{16}\ln^2{2} - \frac12\sum_{k=1}^\infty \frac{(-1)^k H_k}{(2k+1)^2}.
\end{align*}
Nevertheless, the irrationality or transcendence of the constant $\Im(\li_3(1+\img))$ still remains unknown. Readers are encouraged to delve deeper into the properties of $\Im(\li_3(1+\img))$ or explore improved approximation schemes for this intriguing constant.

\section*{Acknowledgment}
I would like to express my gratitude to the Spirit of Ramanujan (SOR) STEM Talent Initiative, directed by Ken Ono, for providing me with the computational tools that were instrumental in verifying and validating my results. I would also like to sincerely thank Cornel I. V\u{a}lean for engaging in valuable discussions with me. Furthermore, I extend my sincere thanks to Ken Ono for his unwavering encouragement throughout this journey.
\section*{Funding}
The author did not receive funding from any organization for the submitted work.
\bibliographystyle{unsrt}
\bibliography{newresearch}
\end{document}